\documentclass{amsart}

\usepackage[all]{xy}

\CompileMatrices

\begin{document}

\newcommand{\nc}{\newcommand}

\nc{\oppA}{A^{\sharp}}
\nc{\oppAo}{A^{\sharp}_0}
\renewcommand{\k}{{\bf k}}
\nc{\e}{{\varepsilon}}
\nc{\ke}{{\bf k}_\e}
\nc{\kb}{{\bf k}_\beta}
\nc{\ra}{\rightarrow}
\nc{\Alm}{A-{\rm mod}}
\nc{\Arm}{{\rm mod}-A}
\nc{\Almb}{(A-{\rm mod})_0}
\nc{\Armb}{({\rm mod}-A)_0}
\nc{\Aormb}{({\rm mod}-A_0)_0}
\nc{\HA}{{\rm Hom}_A}
\nc{\hA}{{\rm hom}_A}
\nc{\hoppA}{{\rm hom}_{\oppA}}
\nc{\CAl}{{\rm Kom}(A)}
\nc{\KAl}{{K}(A)}
\nc{\DAl}{{D}(A)}
\nc{\CAlb}{{\rm Kom}(A)_0}
\nc{\KAlb}{{K}(A)_0}
\nc{\DAlb}{{D}(A)_0}
\nc{\KoppAlb}{{K}(\oppA)_0}
\nc{\DoppAlb}{{D}(\oppA)_0}
\nc{\HDAl}{{\rm Hom}_{\DAl}^{{\gr}}}
\nc{\HKAl}{{\rm Hom}_{\KAl}^{{\gr}}}
\nc{\hDAl}{{\rm hom}_{\DAlb}^{{\gr}}}
\nc{\hKAl}{{\rm hom}_{\KAlb}^{{\gr}}}
\nc{\hDoppAl}{{\rm hom}_{\DoppAlb}^{{\gr}}}
\nc{\hKoppAl}{{\rm hom}_{\KoppAlb}^{{\gr}}}
\nc{\HAd}{\HA^{{\gr}}}
\nc{\hAd}{\hA^{{\gr}}}
\nc{\Hk}{{\rm Hk}^{\gr}(A,A_0,\e )}
\nc{\Hks}{{\rm Hk}^{\frac{\infty}{2}+\gr}(A,A_0,\e )}
\nc{\Hkso}{{\rm Hk}^{\frac{\infty}{2}+\gr}}
\nc{\pr}{\noindent{\em Proof. }}
\nc{\g}{\mathfrak g}
\nc{\ag}{\widehat{\mathfrak g}}
\nc{\an}{\tilde{\mathfrak n}}
\nc{\Uag}{U(\widehat{\mathfrak g})}
\nc{\cUag}{\widehat U(\widehat{\mathfrak g})}
\nc{\Uan}{U(\tilde{\mathfrak n})}
\nc{\Uanp}{U(\tilde{\mathfrak n}_+)}
\nc{\Uanm}{U(\tilde{\mathfrak n}_-)}
\nc{\Uagp}{U(\widehat{\mathfrak g}_+)}
\nc{\Uagm}{U(\widehat{\mathfrak g}_-)}
\nc{\oppUag}{U(\widehat{\mathfrak g})^\sharp}
\nc{\oppUan}{U(\tilde{\mathfrak n})^\sharp}
\nc{\oppg}{{\g}^\sharp}
\nc{\n}{\mathfrak n}
\nc{\h}{\mathfrak h}
\renewcommand{\b}{\mathfrak b}
\nc{\Ug}{U(\g)}
\nc{\Ugo}{U(\h)}
\nc{\oppUgo}{U(\h)^\sharp}
\nc{\oppUg}{U(\g)^\sharp}
\nc{\Uh}{U(\h)}
\nc{\Un}{U(\n)}
\nc{\Ub}{U(\b)}
\nc{\Ubo}{U(\h_{\b})}
\nc{\Uno}{U(\h_{\n})}
\nc{\Ugrmb}{({\rm mod}-\Ug)_0}
\nc{\KoppUaglb}{{K}(\oppUag_k)_0}
\nc{\DoppUaglb}{{D}(\oppUag_k)_0}
\nc{\CoppUglb}{{\rm Kom}(\oppUg)_0}
\nc{\hoppUag}{{\rm hom}_{\oppUag}}
\nc{\hDoppUagl}{{\rm hom}_{\DoppUaglb}^{{\gr}}}
\nc{\hKoppUagl}{{\rm hom}_{\KoppUaglb}^{{\gr}}}
\nc{\hoppUagd}{\hoppUag^{\gr}}
\nc{\gr}{\bullet} 
\nc{\oppAlmb}{(\oppA-{\rm mod})_0}
\nc{\oppArmb}{({\rm mod}-\oppA )_0}
\nc{\oppAlm}{\oppA-{\rm mod}}
\nc{\spr}{\otimes^N_B}
\nc{\spro}{\otimes^{N_0}_{B_0}}
\nc{\stor}{{\rm Tor}_A^{\frac{\infty}{2}+\gr}}
\nc{\storg}{{\rm Tor}_{\Ug}^{\frac{\infty}{2}+\gr}}
\nc{\storo}{{\rm Tor}_{A_0}^{\frac{\infty}{2}+\gr}}
\nc{\storn}{{\rm Tor}_A^{\frac{\infty}{2}+n}}
\renewcommand{\Bar}{{\rm Bar}^{\gr}}
\nc{\Barn}{{\rm Bar}^{-n}}
\nc{\Baro}{{\rm Bar}}
\nc{\tilBarn}{\widetilde{\rm Bar}^{-n}}
\nc{\tilBar}{\widetilde{\rm Bar}^{\gr}}
\nc{\linBarn}{\overline{\rm Bar}^{-n}}
\nc{\linBar}{\overline{\rm Bar}^{\gr}}
\nc{\sBar}{{\rm Bar}^{\frac{\infty}{2}+\gr}}
\nc{\sBarn}{{\rm Bar}^{\frac{\infty}{2}+n}}  
\nc{\sBaropp}{{\rm Bar}^{\frac{\infty}{2}+\gr}_{\sharp}}
\nc{\sBarnopp}{{\rm Bar}^{\frac{\infty}{2}+n}_{\sharp}}
\renewcommand{\sl}{\Lambda^{\frac{\infty}{2}+\gr}(\g)}  
\nc{\slo}{\Lambda^{\frac{\infty}{2}+\gr}(\h)}
\nc{\sloo}{\Lambda^{\frac{\infty}{2}+\gr}}
\nc{\spran}{\otimes^{U(\an_+)}_{U(\an_-)}} 

\newtheorem{theorem}{Theorem}{}
\newtheorem{lemma}[theorem]{Lemma}{}
\newtheorem{corollary}[theorem]{Corollary}{}
\newtheorem{conjecture}[theorem]{Conjecture}{}
\newtheorem{proposition}[theorem]{Proposition}{}
\newtheorem{axiom}{Axiom}{}
\newtheorem{remark}{Remark}{}
\newtheorem{example}{Example}{}
\newtheorem{exercise}{Exercise}{}
\newtheorem{definition}{Definition}{}

\renewcommand{\theproposition}{\thesubsection.\arabic{proposition}}

\renewcommand{\thelemma}{\thesubsection.\arabic{lemma}}

\renewcommand{\thecorollary}{\thesubsection.\arabic{corollary}}

\renewcommand{\theremark}{}

\renewcommand{\thedefinition}{\arabic{definition}}

\renewcommand{\thetheorem}{\thesubsection.\arabic{theorem}}

\title{Semi--infinite cohomology and Hecke algebras}

\author[A. Sevostyanov]{A. Sevostyanov} 
\address{Max-Planck-Institute f\"{u}r Mathematik \\
Box 7280  \\ \mbox{D-53072} Bonn \\ Germany
}
\email{sevastia@mpim-bonn.mpg.de}


\begin{abstract}
This paper provides a homological algebraic foundation for generalizations of classical Hecke algebras introduced in \cite{S}. These new Hecke algebras are associated to triples of the form $(A,A_0,\e)$, where $A$ is an associative algebra over a field $\k$ containing subalgebra $A_0$ with augmentation $\e:A_0\ra \k$. 

These algebras are connected with cohomology of associative algebras in the sense that for every left $A$--module $V$ and right $A$--module $W$ the Hecke algebra associated to triple $(A,A_0,\e)$  naturally acts in the $A_0$--cohomology and $A_0$--homology spaces of $V$ and $W$, respectively. 

We also introduce the semi--infinite cohomology functor for associative algebras and define modifications of Hecke algebras acting in semi--infinite cohomology spaces. We call these algebras semi--infinite Hecke algebras.

As an example we realize the W--algebra $W_k(\g)$ associated to a complex semisimple Lie algebra $\g$ as a semi--infinite Hecke algebra.  Using this realization we explicitly calculate the algebra $W_k(\g)$ avoiding the bosonization technique used in \cite{FF}.
\end{abstract}

\maketitle

\tableofcontents


\section*{Introduction}

\renewcommand{\theequation}{\arabic{equation}}

\setcounter{equation}{0}

Let $G$ be a Chevalley group over a finite field, $B$ a Borel subgroup of $G$, ${\bf 1}_B$ the trivial complex representation of $B$. Denote by $G\otimes_B{\bf 1}_B$ the induced representation of the group $G$. The algebra
\begin{equation}\label{clhk}
{\rm Hom}_G(G\otimes_B{\bf 1}_B,G\otimes_B{\bf 1}_B)
\end{equation}
is called the Hecke algebra of the triple $(G,B,{\bf 1}_B)$ (see, for instance, \cite{I}). 

At present Hecke algebras play an important role in various fields of mathematics, e.g. in representation theory of Chevalley groups (see \cite{Bo,C} and references there). Among of the other applications one should mention the Kazhdan--Lusztig polynomials \cite{KL}. 

It turns out that algebras of a similar type appear in the study of algebraic objects of different nature. The purpose of this paper is to investigate general properties of these algebras. 

First we shall define an abstract version of the classical Hecke algebra (\ref{clhk}) in the following general situation. 
Let $A$ be an associative algebra over a field $\k$,
$A_0\subset A$ a subalgebra with augmentation $\e : A_0 \ra \k$. We denote this one--dimensional $A_0$--module by $\ke$.
Let $X^{\gr}$ be a projective resolution of the left $A_0$--module $\ke$. Consider the complex $A\otimes_{A_0}X^{\gr}$ of left $A$--modules, where $A_0$ acts on $A$ by right multiplication and the structure of a left $A$--module on $A\otimes_{A_0}X^{\gr}$ is induced by the left regular action of $A$. We call the ${\mathbb Z}$--graded algebra
\begin{equation}\label{heckedefo}
\Hk = \bigoplus_{n\in {\mathbb Z}}{\rm Hom}_{D(A)}(A\otimes_{A_0}X^{\gr},T^n(A\otimes_{A_0}X^{\gr})),
\end{equation}
where $D(A)$ is the derived category of the category of left $A$--modules and $T$ is the grading shift functor, the Hecke algebra of the triple $(A,A_0,\e)$.

Note that if $H^{\gr}(A\otimes_{A_0}X^{\gr})={\rm Tor}_{A_0}^{\gr}(A,\ke )=A\otimes_{A_0}\ke$ the zeroth graded component of the algebra $\Hk$ takes the form
\begin{equation}\label{hkoo}
{\rm Hk}^0(A,A_0,\e)=\HA(A\otimes_{A_0}\ke,A\otimes_{A_0}\ke).
\end{equation}
Therefore the Hecke algebra of the triple $(A,A_0,\e)$ is a natural generalization of the classical Hecke algebra (\ref{clhk})\footnote{The situation described here often appears in quantum mechanics and quantum field theory. In physical vocabulary the triple $(A,A_0,\e)$ is called a quantum system with first--class constraints. The algebra (\ref{hkoo}) is the corresponding quantum reduced system (see, for instance, \cite{D} or discussion in \cite{S}, Sect. 6).}. 

Particular examples of algebras of the form (\ref{hkoo}) are the realization of the center of the universal enveloping algebra of a complex semisimple Lie algebra obtained by Kostant in \cite{K} and algebras of invariant differential operators on homogeneous spaces described in \cite{duf1,duf2,Koorw}. In Section \ref{ex} we discuss these examples in detail. 

The definition (\ref{heckedefo}) is motivated by the  quantum BRST reduction procedure. In \cite{S} the algebra $\Hk$ was defined as the cohomology of a differential graded algebra that is a generalization of the quantum BRST complex proposed in \cite{KSt} (see \cite{S}, Sect. 5). The quantum BRST complex was constructed in \cite{KSt} for triples of the type $(A,U(\g_0),\e)$, where $U(\g_0)$ is the universal enveloping algebra of a finite--dimensional Lie algebra, and $\e$ is the trivial representation of $U(\g_0)$. Remarkably, this complex already appeared in unpublished lectures \cite{duf1} in the situation when $A=\Ug$, where $\g$ is a finite--dimensional Lie algebra containing $\g_0$ as a subalgebra.

Hecke algebras of the form (\ref{heckedefo}) are associated to the usual cohomology of associative algebras in the sense that for every left $A$--module $V$ and right $A$--module $W$ the algebra $\Hk$ naturally acts in the $A_0$--cohomology and $A_0$--homology spaces of $V$ and $W$, $H^{\gr}(A_0,V)$ and $H_{\gr}(A_0,W)$, from the right and from the left, respectively (see Proposition \ref{hact}); here $A_0$ is augmented by $\e$. 

In this paper we also define modifications of Hecke algebras acting in semi--infinite cohomology spaces. The semi--infinite cohomology was first introduced in \cite{F}, for a class of $\mathbb Z$--graded Lie algebras with finite--dimensional graded components, as the cohomology of a standard complex. In \cite{V} (see also \cite{V1,V2}) this definition was explained from the point of view of homological algebra: the semi--infinite cohomology was obtained as a two--sided derived functor of the functor of semivariants which is a mixture of the functor of invariants and of the functor of covariants.

In \cite{Arkh1}--\cite{Arkh3} using a kind of Koszul duality S. Arkhipov tries to generalize the semi--infinite cohomology functor to a class of $\mathbb Z$--graded associative algebras. However papers \cite{Arkh1}--\cite{Arkh3} contain numerous mistakes\footnote{See, for instance the footnote to Lemma \ref{abovebound} in this paper}. As a consequence the definition of the semi--infinite Tor functor given in \cite{Arkh3} is not self--consistent and the  relation of this functor to the semi--infinite cohomology is not clear. 

In this paper we give a correct definition of the semi--infinite Tor functor for associative algebras. Our definition is a direct generalization of the original Voronov's definition: the semi--infinite Tor functor is defined as a derived functor of the functor of semiproduct, which is, in turn, a generalization of the functor of semivariants to the case of associative algebras. We also derive some new properties of the semi--infinite cohomology functor (see Theorem \ref{threetor}). The proofs of these properties are quite difficult technically, and we moved them to the Appendix.

In the construction of the semi--infinite Tor functor we use the notion of the semiregular bimodule that plays the role of the left and right regular representations in the semi--infinite cohomology theory. The semiregular bimodule was introduced in \cite{V} (see also \cite{Sor}) in the Lie algebra case. In this paper we also use some correct results of \cite{Arkh1} on semiregular bimodules for associative algebras.

Finally we define a modification of Hecke algebras (\ref{heckedefo}) associated to semi--infinite cohomology. We call these new algebras the semi--infinite Hecke algebras.  The example discussed in Section \ref{Hklie} shows that the semi--infinite Hecke algebras, as well as the semi--infinite cohomology, are adapted for the study of infinite--dimensional objects. Our main result in Section \ref{Hklie}  is a realization of the W--algebra $W_k(\g)$ associated to a complex semisimple Lie algebra $\g$ (see \cite{FF} for the definition of this algebra) as a semi--infinite Hecke algebra (see Proposition \ref{w}).  Using this realization we explicitly calculate the algebra $W_k(\g)$ (see Theorem \ref{mw}) avoiding the bosonization technique used in \cite{FF}. The description of the algebra $W_k(\g)$ obtained in Theorem \ref{mw} is similar to the Kostant's realization of the center $Z(\Ug)$ of the universal enveloping algebra $\Ug$.
\vskip 0.3cm
\noindent
{\bf Acknowledgments.}
I am very grateful to M. Duflo who sent me unpublished lectures \cite{duf1}, Yu. Manin for useful discussion, A. Vishik for important remarks on paper \cite{S} and S. Yagunov who found a time for numerous discussions of complicated technical details of this paper.

\renewcommand{\theequation}{\thesubsection.\arabic{equation}}


\section{Hecke algebras}

In this section, using the language of derived categories, we give a definition of Hecke algebras. We also show that this definition is equivalent to the elementary definition of Hecke algebras proposed in \cite{S}. However our new definition will be useful for generalization of the notion of Hecke algebras to semi--infinite cohomology discussed in Section \ref{shecke}. Using this definition we also obtain elementary proofs of the properties of Hecke algebras derived in \cite{S}.


\subsection{Notation and recollections}\label{notat}

\setcounter{equation}{0}
\setcounter{theorem}{0}

Let ${\mathcal A}$ be an abelian category.
In this section we recall, following \cite{GM}, main facts about homotopy and derived categories associated to ${\mathcal A}$. These facts will be used throughout of this paper.

Let ${{\rm Kom}({\mathcal A})}$ be the category of complexes over ${\mathcal A}$, ${{K}({\mathcal A})}$ the corresponding homotopy category. The category ${{K}({\mathcal A})}$ has the same objects as ${{\rm Kom}({\mathcal A})}$, morphisms of ${{K}({\mathcal A})}$ being morphisms of complexes modulo homotopic equivalence. We denote by ${{D}({\mathcal A})}$ the derived category of the category ${\mathcal A}$. ${{D}({\mathcal A})}$ is the localization of the homotopy category ${{K}({\mathcal A})}$ by the class of quasi--isomorphisms (see \cite{GM}, Ch. III).

Any object $X$ of the category $\mathcal A$ may be considered as a complex $\ldots\ra 0\ra X\ra 0\ra \ldots$ (with $X$ at the 0-th place). Such complexes are called 0--complexes. The functor ${\mathcal A}\ra K({\mathcal A})$ sending every object of $\mathcal A$ to the corresponding 0--complex is fully faithful. Using this functor we shall always identify $\mathcal A$ with the subcategory of 0--complexes in $K({\mathcal A})$.

A complex $X^{\gr}$ is called an $H^0$--complex if $H^i(X^{\gr})=0$ for $i\neq 0$. Such complexes form a full subcategory in $D({\mathcal A})$. The following proposition shows that $\mathcal A$ may be regarded not only as a subcategory in $K({\mathcal A})$ but also as a subcategory in $D({\mathcal A})$.

\begin{proposition}{\bf (\cite{GM}, Proposition III.5.2.)}\label{ocompl}
The functor ${\mathcal A}\ra D({\mathcal A})$ sending every object of $\mathcal A$ to the corresponding 0--complex yields an equivalence of $\mathcal A$ with the full subcategory of $D({\mathcal A})$ formed by $H^0$--complexes.
\end{proposition}

We shall use the graded $\rm Hom$ in the category ${{D}({\mathcal A})}$ introduced by 
\begin{equation}
{{\rm Hom}_{{D}({\mathcal A})}^{{\gr}}}(X^{\gr},Y^{\gr})=\bigoplus_{n\in {\mathbb Z}}{\rm Hom}_{{{D}({\mathcal A})}}(X^{\gr},Y[n]^{\gr}),
\end{equation}
where the complex $Y[n]^{\gr}$ is defined by
$$
Y[n]^k=Y^{k+n},~d_{Y[n]^{\gr}}=(-1)^nd_{Y^{\gr}}.
$$

Similarly we define
\begin{equation}
{{\rm Hom}_{{K}({\mathcal A})}^{{\gr}}}(X^{\gr},Y^{\gr})=\bigoplus_{n\in {\mathbb Z}}{\rm Hom}_{{{K}({\mathcal A})}}(X^{\gr},Y[n]^{\gr}).
\end{equation}

Recall that the space ${{\rm Hom}_{{K}({\mathcal A})}^{{\gr}}}(X^{\gr},Y^{\gr})$ may be calculated as follows (see \cite{GM}, III.6.14). Consider a complex ${\rm Hom}_{\mathcal A}^{\gr}(X^{\gr},Y^{\gr})$,
$$
\begin{array}{l}
{\rm Hom}_{\mathcal A}^{\gr}(X^{\gr},Y^{\gr})=\bigoplus_{n\in {\mathbb Z}}{\rm Hom}_{\mathcal A}^n(X^{\gr},Y^{\gr}), \\
 \\
{\rm Hom}_{\mathcal A}^n(X^{\gr},Y^{\gr})=\prod_{p\in {\mathbb Z}}{\rm Hom}_{\mathcal A}(X^p,Y^{p+n}) \\
\end{array}
$$
with the differential given by
$$
{\bf d}f=d_{Y^{\gr}}\circ f-(-1)^{n}f\circ d_{X^{\gr}},~~f\in {\rm Hom}_{\mathcal A}^n(X^{\gr},Y^{\gr}).
$$
Then
\begin{equation}\label{homK}
{{\rm Hom}_{{K}({\mathcal A})}^{{\gr}}}(X^{{\gr}},Y^{{\gr}})=H^{\gr}({\rm Hom}_{\mathcal A}^{\gr}(X^{\gr},Y^{\gr})).
\end{equation}

The main property of derived categories is that 
sometimes they  may be realized as homotopy categories. For instance, let
$D^+({\mathcal A})$ ($D^-({\mathcal A})$) be the full subcategory in ${{D}({\mathcal A})}$ whose objects are complexes bounded from below (above). Let ${\mathcal I}({\mathcal A})$ (${\mathcal P}({\mathcal A})$) be the full subcategory in ${\mathcal A}$ formed by injective (projective) objects, ${\rm Kom}^+({\mathcal I}({\mathcal A}))$ (${\rm Kom}^-({\mathcal P}({\mathcal A}))$) the category of complexes bounded from below (above) over this abelian category, $K^+({\mathcal I}({\mathcal A}))$ ($K^-({\mathcal P}({\mathcal A}))$) the corresponding homotopy category.

\begin{proposition}{\bf (\cite{GM}, Theorem III.5.21)}\label{equiv}
Suppose that the category $\mathcal A$ has enough injective and projective objects, i.e. for every object $X\in {\rm Ob}~{\mathcal A}$ there exist an injection into an injective object, $X\ra I,~I\in {\rm Ob}~{\mathcal I}({\mathcal A})$, and a surjection $P\ra X$ from a projective object $P\in {\rm Ob}~{\mathcal P}({\mathcal A})$ onto $X$. Then   
the functor of localization by the class of quasi--isomorphisms is an equivalence of categories:
$$
\begin{array}{l}
K^+({\mathcal I}({\mathcal A}))\ra D^+({\mathcal A}), \\
\\
K^-({\mathcal P}({\mathcal A}))\ra D^-({\mathcal A}).
\end{array}
$$

Moreover, let $X^{\gr},~Y^{\gr}$ be two objects of the category ${{K}({\mathcal A})}$ such that  either $Y^{\gr}\in {\rm Ob}~K^+({\mathcal I}({\mathcal A}))$ or $X^{\gr}\in {\rm Ob}~K^-({\mathcal P}({\mathcal A}))$. Then the natural homomorphism
$$
{{\rm Hom}_{{K}({\mathcal A})}^{{\gr}}}(X^{\gr},Y^{\gr})\ra {{\rm Hom}_{{D}({\mathcal A})}^{{\gr}}}(X^{\gr},Y^{\gr})
$$
is an isomorphism.
\end{proposition}


\subsection{Hecke algebras: definition and main properties}\label{hecke}

\setcounter{equation}{0}
\setcounter{theorem}{0}

In this section we introduce the notion of Hecke algebras which is the main object of our study in this paper. 

Let $A$ be an associative algebra over a field $\k$.
Let $\Alm$ be the category of left $A$--modules. We denote by $\HA (\cdot , \cdot)$ the set of morphisms between two objects of this category. We also write $K(A)$, $D(A)$ and  ${\HA^{{\gr}}}$ instead of $K(\Alm)$, $D(\Alm)$ and
${\rm Hom}_{\Alm}^{\gr}$, respectively.

Suppose that $A$ contains a subalgebra $A_0$ with augmentation $\e : A_0 \ra \k$. We denote this one--dimensional $A_0$-- module by $\ke$.

\begin{definition}
Let $X^{\gr}$ be a projective resolution of the left $A_0$--module $\ke$. Consider the complex $A\otimes_{A_0}X^{\gr}$ of left $A$--modules, where $A_0$ acts on $A$ by right multiplication and the structure of a left $A$--module on $A\otimes_{A_0}X^{\gr}$ is induced by the left regular action of $A$. The graded algebra
\begin{equation}\label{heckedef}
\Hk = \HDAl(A\otimes_{A_0}X^{\gr},A\otimes_{A_0}X^{\gr})
\end{equation}
is called the Hecke algebra of the triple $(A,A_0,\e)$.
\end{definition}
Clearly, the algebra defined by (\ref{heckedef}) does not depend on the choice of the resolution $X^{\gr}$.

The following vanishing property allows to explicitly calculate Hecke algebras in many particular situations. 

\begin{proposition}{\bf (\cite{S}, Theorem 7)}\label{vanish}

Assume that $H^{\gr}(A\otimes_{A_0}X^{\gr})={\rm Tor}_{A_0}^{\gr}(A,\ke )=A\otimes_{A_0}\ke$. Then
$$
\Hk=\HDAl(A\otimes_{A_0}\ke,A\otimes_{A_0}\ke).
$$

In particular,
$$
{\rm Hk}^0(A,A_0,\e)=\HA(A\otimes_{A_0}\ke,A\otimes_{A_0}\ke).
$$
\end{proposition}
\begin{remark}
The condition ${\rm Tor}_{A_0}^{\gr}(A,\ke )=A\otimes_{A_0}\ke$ is satisfied, for instance, if $A$ is projective as a right $A_0$--module.
\end{remark}
\pr The condition $H^{\gr}(A\otimes_{A_0}X^{\gr})=A\otimes_{A_0}\ke$ implies that the natural map $A\otimes_{A_0}X^{\gr}\ra A\otimes_{A_0}\ke$ is a quasi--isomorphism, i.e. this is an isomorphism in the category $\DAl$. Therefore in (\ref{heckedef}) we can replace the complex $A\otimes_{A_0}X^{\gr}$ with the 0--complex that corresponds to the left $A$--module $A\otimes_{A_0}\ke$. This completes the proof.
\qed

Another important property of Hecke algebras is that for every left $A$--module $V$ and right $A$--module $W$ the algebra $\Hk$ naturally acts in the $A_0$--cohomology and $A_0$--homology spaces of $V$ and $W$, $H^{\gr}(A_0,V)$ and $H_{\gr}(A_0,W)$, $A_0$ being augmented by $\e$.

First we shall define an action of the  Hecke algebra $\Hk$ in the  cohomology space $H^{\gr}(A_0,V)$. 
Recall that the cohomology module $H^{\gr}(A_0,V)$ may be defined as the cohomology space of the complex ${\rm Hom}_{A_0}^{\gr}(X^{\gr},V)$, where $X^{\gr}$ is a projective resolution of the left $A_0$--module $\ke$. The complex ${\rm Hom}_{A_0}^{\gr}(X^{\gr},V)$ is canonically isomorphic to the complex $\HAd(A\otimes_{A_0}X^{\gr},V)$, and hence
$$
H^{\gr}(A_0,V)=H^{\gr}(\HAd(A\otimes_{A_0}X^{\gr},V)).
$$

Using (\ref{homK}) we can also represent this expression for the cohomology space $H^{\gr}(A_0,V)$ in the following equivalent form:
\begin{equation}\label{cohomol}
H^{\gr}(A_0,V)=\HKAl(A\otimes_{A_0}X^{\gr},V).
\end{equation}
 
Remark that since $A$ is $A$--projective as a left $A$--module, and $X^{\gr}$ is a bounded from above complex of $A_0$--projective modules, $A\otimes_{A_0}X^{\gr}$ is a bounded from above complex of $A$--projective modules, i.e $A\otimes_{A_0}X^{\gr}$ is an object of the category $K^-({\mathcal P}(\Alm))$. Therefore by Proposition \ref{equiv} we have:
$$
H^{\gr}(A_0,V)=\HDAl(A\otimes_{A_0}X^{\gr},V).
$$

Finally observe that the Hecke algebra $\Hk$ naturally acts on the r.h.s. of the last equality by compositions of homomorphisms.

Now we turn to the definition of an action of Hecke algebras on homology spaces.
The homology space $H_{\gr}(A_0,W)$ may be calculated as follows:
\begin{equation}\label{homol}
H_{\gr}(A_0,W)=H^{\gr}(W\otimes_{A_0}X^{\gr})=H^{\gr}(W\otimes_AA\otimes_{A_0}X^{\gr}).
\end{equation}

Since $A\otimes_{A_0}X^{\gr}$ is an object of the category $K^-({\mathcal P}(\Alm))$, the last expression may be rewritten using 
the definition of the left derived functor $\otimes_A^L$ of the functor $\otimes_A$ (see \cite{GM}, III.6.15 and III.7, Ex. 6) as follows:
$$
H_{\gr}(A_0,W)=H^{\gr}(W\otimes_A^L(A\otimes_{A_0}X^{\gr})).
$$

The algebra $\Hk$ naturally acts on the r.h.s. of the last expression by applying of endomorphisms to the second argument $A\otimes_{A_0}X^{\gr}$ of the derived functor $\otimes_A^L$.

Clearly, the actions defined above respect the gradings of $\Hk,~H^{\gr}(A_0,V)$ and $H_{\gr}(A_0,W)$. Thus we have proved the following statement.

\begin{proposition}{\bf (\cite{S}, Theorems 5 and 6)}\label{hact}
For every left $A$--module $V$ and right $A$--module $W$ the algebra $\Hk$ naturally acts in the $A_0$--cohomology and $A_0$--homology spaces of $V$ and $W$, $H^{\gr}(A_0,V)$ and $H_{\gr}(A_0,W)$, from the right and from the left, respectively. These actions respect the gradings of $\Hk$,\\ $H^{\gr}(A_0,V)$ and $H_{\gr}(A_0,W)$, i.e.
$$
\begin{array}{l}
H^n(A_0,V)\times {\rm Hk}^m(A,A_0,\e )\ra H^{n+m}(A_0,V), \\
\\
{\rm Hk}^m(A,A_0,\e )\times H_n(A_0,W)\ra H_{n-m}(A_0,W).
\end{array}
$$
\end{proposition}

In conclusion we note that using Proposition \ref{equiv}, formula (\ref{homK}) and the fact that $A\otimes_{A_0}X^{\gr}$ is an object of the category $K^-({\mathcal P}(\Alm))$  the definition of the Hecke algebra may be rewritten as follows:
\begin{equation}\label{homdef}
\Hk=H^{\gr}(\HAd(A\otimes_{A_0}X^{\gr},A\otimes_{A_0}X^{\gr})).
\end{equation}
It is the definition of Hecke algebras that first appeared in \cite{S}.


\subsection{Examples of Hecke algebras}\label{ex}

\setcounter{equation}{0}
\setcounter{theorem}{0}

In the Introduction we already observed that the Hecke algebras introduced in this paper are generalizations of the classical Hecke algebras associated to Chevalley groups over finite fields.
In this section we discuss in detail the other important examples of Hecke algebras mentioned in the Introduction.
\vskip 0.3cm
\noindent
{\bf Example 1. The center of the universal enveloping algebra of a complex semisimple Lie algebra}
\vskip 0.3cm
Let $\g$ be a complex semisimple Lie algebra, $\n\subset \g$ a maximal nilpotent subalgebra, $\Ug$ and $\Un$ the universal enveloping algebras of $\g$ and $\n$, respectively. Denote by $X_i,~i=1,\ldots,{\rm rank}~\g$ simple root vectors in $\n$. 

Let $\chi:\n\ra {\mathbb C}$ be a character of $\n$. We denote the corresponding one--dimensional $\Un$--module by ${\mathbb C}_\chi$.  Since $\n=\sum_{i=1}^{{\rm rank}\g}{\mathbb C}X_i\oplus [\n,\n]$ the character $\chi$ is completely determined by the constants $\chi(X_i),~i=1,\ldots,{\rm rank}~\g$. Such a character is called non--singular if all these constants are not equal to zero.

\begin{proposition}{\bf (\cite{K}, Theorem 2.4.2)}
Suppose that $\chi:\n\ra {\mathbb C}$ is a non--singular character of a maximal nilpotent subalgebra ${\n}\subset \g$. Then the algebra ${\rm End}_{\Ug}(\Ug\otimes_{\Un}{\mathbb C}_\chi)^{opp}$ is canonically isomorphic to the center $Z(\Ug)$ of the universal enveloping algebra $\Ug$,
$$
{\rm End}_{\Ug}(\Ug\otimes_{\Un}{\mathbb C}_\chi)^{opp}= Z(\Ug).
$$
\end{proposition}

Now consider the Hecke algebra of the triple $(\Ug,\Un,\chi)$. 
Since $\Ug$ is projective as a right $\Un$--module (see \cite{carteil}, Ch. XIII, Proposition 4.1) the conditions of Proposition \ref{vanish} are satisfied, and the algebra ${\rm Hk}^0(\Ug,\Un,\chi)$ is isomorphic to ${\rm End}_{\Ug}(\Ug\otimes_{\Un}{\mathbb C}_\chi)$,
$$
{\rm End}_{\Ug}(\Ug\otimes_{\Un}{\mathbb C}_\chi)= {\rm Hk}^0(\Ug,\Un,\chi).
$$

Thus the center $Z(\Ug)$ is realized as the zeroth graded component of the Hecke algebra of the triple $(\Ug,\Un,\chi)$,
$$
{\rm Hk}^0(\Ug,\Un,\chi)^{opp}= Z(\Ug).
$$
This result, as well as a similar realization of the center of a quantum group, was obtained in \cite{S1}.
\vskip 0.3cm
\noindent
{\bf Example 2. Algebras of invariant differential operators on homogeneous spaces}
\vskip 0.3cm
Let $\g$ be a finite--dimensional Lie algebra over $\mathbb R$, $\h\subset \g$ a subalgebra in $\g$, $\Ug$ and $\Uh$ the universal enveloping algebras of $\g$ and $\h$, respectively. Let $\chi:\h\ra {\mathbb R}$ be a character of $\h$. We denote the corresponding one--dimensional $\Uh$--module by ${\mathbb R}_{\chi}$. 

Let $G$ be a connected group with Lie algebra $\g$, $H\subset G$ a closed connected subgroup with Lie algebra $\h$. Let $L_\chi=G\times_H {\mathbb R}_{\chi}$ be the line bundle on the homogeneous space $G/H$ associated to ${\mathbb R}_{\chi}$.
The algebra $D_\chi={\rm End}_{\Ug}(\Ug\otimes_{\Uh}{\mathbb R}_\chi)$ is isomorphic to the algebra of $G$--invariant differential operators acting on the space of smooth sections of $L_\chi$ (see \cite{Koorw, duf1, duf2}).

Now consider the Hecke algebra of the triple $(\Ug,\Uh,\chi)$.
Since $\Ug$ is projective as a right $\Uh$--module (see \cite{carteil}, Ch. XIII, Proposition 4.1) the conditions of Proposition \ref{vanish} are satisfied, and  we have an algebraic isomorphism:
$$
{\rm End}_{\Ug}(\Ug\otimes_{\Uh}{\mathbb R}_\chi)={\rm Hk}^0(\Ug,\Uh,\chi),
$$
i.e. the algebra $D_\chi$ is realized as the zeroth graded component of the Hecke algebra of the triple $(\Ug,\Uh,\chi)$,
$$
D_\chi= {\rm Hk}^0(\Ug,\Uh,\chi).
$$


\section{Semi-infinite cohomology}

In this section, using results of semi--infinite homological algebra (see \cite{V}), we shall define the semi--infinite $\rm Tor$ functor for a class of graded associative algebras. The properties of this functor turn out to be quite similar to those of the usual $\rm Tor$ functor. The semi--infinite Tor functor is a generalization of the functor of semi--infinite cohomology introduced in \cite{F,V} for a class of graded Lie algebras.


\subsection{Notation and conventions}\label{setup}

\setcounter{equation}{0}
\setcounter{theorem}{0}

We shall define the semi--infinite $\rm Tor$ functor for a class of $\mathbb Z$--graded associative algebras over a field $\k$. Let $A$ be such an algebra,
$$
A=\bigoplus_{n\in {\mathbb Z}}A_n.
$$
The category of left (right) $\mathbb Z$--graded $A$--modules with morphisms being homomorphisms of $A$--modules preserving gradings is denoted by $\Alm$ ($\Arm$). For both of these categories the set of morphisms between two objects is denoted by $\HA(\cdot,\cdot)$. For $M,M^\prime\in {\rm Ob}~\Alm~({\rm Ob}~\Arm)$ we shall also frequently use the space of homomorphisms of all possible degrees with respect to the gradings on $M$ and $M^\prime$ introduced by
$$
\hA(M,M^\prime)=\bigoplus_{n\in {\mathbb Z}}\HA(M,M^\prime\langle n\rangle ),
$$
where the module $M^\prime\langle n\rangle $ is obtained from $M^\prime$ by grading shift as follows:
$$
M^\prime\langle n\rangle _k=M^\prime_{k+n}.
$$

In this paper we shall mainly deal with the full subcategory of $\Alm$ ($\Arm$) whose objects are modules $M\in {\rm Ob}~\Alm~({\rm Ob}~\Arm)$ such that their gradings  are bounded from above, i.e. 
$$
M=\bigoplus_{n\leq K(M)}M_n,~K(M)\in {\mathbb Z}.
$$
This subcategory is denoted by $\Almb$ ($\Armb$).
We also denote by ${\rm Vect}_{\k}$ the category of $\mathbb Z$--graded vector spaces over $\k$. 

All tensor products of graded $A$--modules and graded vector spaces will be understood in the graded sense.

The following simple lemma, that is a direct consequence of definitions, will be often used in this paper.
\begin{lemma}\label{tenshom}
Let $M$ and $M^\prime$ be two objects of the category ${\rm Vect}_{\k}$ such that $M=\bigoplus_{n\leq K}M_n,~K\in {\mathbb Z}$, $M^\prime =\bigoplus_{n\geq L}M^\prime_n,~L\in {\mathbb Z}$, and for every $n$ ${\rm dim}~M^\prime_n<\infty$. Then
$$
{\rm hom}_{\k}(M^\prime,M)={M^\prime}^*\otimes M,~\mbox{where }{M^\prime}^*={\rm hom}_{\k}(M^\prime,\k).
$$ 
\end{lemma}

We denote by
$\CAlb$, $\KAlb$ and $\DAlb$ the category of complexes over $\Almb$, the corresponding homotopy and derived category, respectively.
We shall use the double graded $\rm Hom$ in the category $\DAlb$ introduced by 
\begin{equation}
\hDAl(X^{\gr},Y^{\gr})=\bigoplus_{m,n\in {\mathbb Z}}{\rm Hom}_{\DAlb}(X^{\gr},Y[n]\langle m\rangle ^{\gr} ),
\end{equation}
where the complex $Y[n]\langle m\rangle ^{\gr}$ is defined by
$$
Y[n]\langle m\rangle _l^k=Y^{k+n}_{m+l},~d_{Y[n]\langle m\rangle ^{\gr}}=(-1)^nd_{Y^{\gr}}.
$$

Similarly we define 
\begin{equation}
\hKAl(X^{\gr},Y^{\gr})=\bigoplus_{m,n\in {\mathbb Z}}{\rm Hom}_{\KAlb}(X^{\gr},Y[n]\langle m\rangle ^{\gr}).
\end{equation}

Using formula (\ref{homK}) the space $\hKAl(X^{\gr},Y^{\gr})$ may be calculated as follows. Consider a complex $\hAd(X^{\gr},Y^{\gr})$,
$$
\begin{array}{l}
\hAd(X^{\gr},Y^{\gr})=\bigoplus_{n\in {\mathbb Z}}\hA^n(X^{\gr},Y^{\gr}), \\
 \\
\hA^n(X^{\gr},Y^{\gr})=\prod_{p\in {\mathbb Z}}\hA(X^p,Y^{p+n}) \\
\end{array}
$$
with the differential given by
\begin{equation}\label{ddd}
{\bf d}f=d_{Y^{\gr}}\circ f-(-1)^{n}f\circ d_{X^{\gr}},~~f\in \hA^n(X^{\gr},Y^{\gr}).
\end{equation}
Then
\begin{equation}\label{homK1}
\hKAl(X^{{\gr}},Y^{{\gr}})=H^{\gr}(\hAd(X^{\gr},Y^{\gr})).
\end{equation}
 
In order to define the semi--infinite $\rm Tor$ functor we have to impose additional restrictions on the algebra $A$ (see \cite{Arkh1}). Namely, in the rest of this paper we  suppose that $A$  satisfies the following conditions:
{\em
\vskip 0.3cm
\qquad(i) $A$ contains two graded subalgebras $N$ and $B$.
\vskip 0.3cm
\qquad(ii) $N$ is positively graded.
\vskip 0.3cm
\qquad(iii) $N_0=\k$.
\vskip 0.3cm
\qquad(iv) $\dim N_n<\infty$ for any $n\in N$.
}
\vskip 0.3cm
In particular $N$ is naturally augmented.
We denote the augmentation ideal
$\bigoplus_{n>0}N_n$ by $\overline{N}$.
We also denote $\overline{B}=B/\k$.
{\em
\vskip 0.3cm
\qquad(v) $B$ is negatively graded.
\vskip 0.3cm
\qquad(vi) The multiplication in $A$ defines isomorphisms of graded vector spaces
\begin{equation}\label{dectr}
 B\otimes N\ra A \mbox{  and } N\otimes B\ra A.
\end{equation}
}

We call the decompositions (\ref{dectr}) the triangular decompositions for the algebra $A$. 
Note that the compositions of the triangular decomposition maps and of their inverse maps 
yield linear mappings
\begin{equation}\label{cmap}
\begin{array}{l}
N\otimes B \ra B\otimes N, \\
\\
B\otimes N \ra N\otimes B.
\end{array}
\end{equation}
{\em
\vskip 0.3cm
\qquad(vii) The the mappings (\ref{cmap}) are continuous in the following sense: for every $m,n \in {\mathbb Z}$ there exist $k_+,k_-\in{\mathbb Z}$ such that
$$ 
N_m\otimes B_n\ra \bigoplus_{k_-\leq k \leq k_+}B_{n-k}\otimes N_{m+k} 
\mbox{ and } 
B_n\otimes N_m\ra \bigoplus_{k_-\leq k \leq k_+}N_{m-k}\otimes B_{n+k}.
$$
}

\subsection{Semiregular bimodule}\label{bimod}

\setcounter{equation}{0}
\setcounter{theorem}{0}

In this section we recall the definition the semiregular bimodule for the algebra $A$. The notion of the semiregular bimodule was introduced by Voronov (see \cite{V}) in the Lie algebra case and generalized in \cite{Arkh1} to the case of graded associative algebras satisfying conditions (i)--(vii) of the previous section. In the semi--infinite cohomology theory this bimodule plays the role of the regular representation. In particular, the semiregular bimodule naturally appears in the definition of the semi--infinite modification of Hecke algebras.

First consider the right graded $N$-module $N^*={\rm hom}_{\k} (N,\k )$, where the action of $N$ on $N^*$ is defined by
$$
(n\cdot f)(n')=f(nn')\mbox{ for any } f\in N^*,~n\in N.
$$
The right $A$--module
$$
S_A=N^*\otimes_NA
$$
is called the right semiregular representation of $A$ (see \cite{V}, Sect 3.2; \cite{Arkh1}, Sect. 3.4). 

Clearly, that $S_A= N^*\otimes B$ as a right $B$-module. The space $S_A= N^*\otimes B$ is non--positively graded, and hence $S_A\in \Armb$.

Now we obtain another realization for the right semiregular representation. 
Consider another right $A$-module
$S_A^{\prime}={\rm hom}_B(A,B)$, where $B$ acts on $A$ and $B$ by right multiplication. The right action of $A$ on the space $S_A^{\prime}$ is given by
$$
(a\cdot f)(a')=f(aa'),~~ f\in {\rm hom}_B(A,B),~a\in A.
$$

\begin{lemma}{\bf (\cite{Arkh1}, Lemma 3.5.1)}\label{phiiso}
Fix a decomposition 
\begin{equation}\label{decomp}
A=N\otimes  B
\end{equation}
provided by the multiplication in $A$. Let $\phi:S_A \ra S_A^{\prime}$ be a map defined by
$$
\phi(f\otimes a)(a')=f((aa')_N)(aa')_B,
$$
where $f\otimes a\in S_A,~a'\in A$ and $aa'=(aa')_N(aa')_B$ is the decomposition (\ref{decomp}) of the element $aa'$. Then $\phi$ is a morphism of right $A$--modules.
\end{lemma}

We shall suppose that the algebra $A$ satisfies the following additional condition:
{\em
\vskip 0.3cm
\qquad(viii) The morphism $\phi:S_A \ra S_A^{\prime}$ constructed in the previous lemma is an isomorphism of right $A$--modules.
\vskip 0.3cm
}
Finally we have two realizations of the right $A$--module $S_A$:
\begin{equation}\label{SA1}
S_A=N^*\otimes_NA,
\end{equation}
and
\begin{equation}\label{SA2}
S_A={\rm hom}_B(A,B).
\end{equation}

Now we define a structure of a left module on $S_A$ commuting with the right semiregular action of $A$. 
First observe that using realizations (\ref{SA1}) and (\ref{SA2}) of the right semiregular representation one can define natural left actions of the algebras $N$ and $B$ on the space $S_A$ induced by the natural left action of $N$ on $N^*$ and the left regular representation of $B$, respectively. Clearly, these actions commute with the right action of the algebra $A$ on $S_A$. Therefore we have natural inclusions of algebras
$$
N\hookrightarrow \hA(S_A,S_A),~~B\hookrightarrow \hA(S_A,S_A).
$$
Denote by $\oppA$ the subalgebra in $\hA(S_A,S_A)$ generated by $N$ and $B$.
    
\begin{proposition}{\bf (\cite{Arkh1}, Corollary 3.3.3, Lemma 3.5.3 and Corollary 3.5.3)}\label{SAopp}
$\oppA$ is a $\mathbb Z$--graded associative algebra satisfying conditions (i)--(vii) of Section \ref{setup}. Moreover, $S_A\in \oppAlmb$ and
\begin{equation}\label{oppSA1}
S_A=\oppA\otimes_NN^*=\qquad \qquad \qquad \qquad
\end{equation}
\begin{equation}\label{oppSA2}
\qquad \qquad \qquad \qquad = {\rm hom}_B(\oppA,B)
\end{equation}
as a left $\oppA$--module.
\end{proposition}

Using Proposition \ref{SAopp} the space $S_A$ is equipped with the structure of  an $\oppA -A$ bimodule. This bimodule is called the  semiregular bimodule associated to the algebra $A$. The left action of the algebra $\oppA$ on the space $S_A$ is called the left semiregular action.


\subsection{Semiproduct}\label{sspr}

\setcounter{equation}{0}
\setcounter{theorem}{0}

In this section we define the functor of semiproduct. This functor is a generalization of the functor of semivariants (see \cite{V}, Sect. 3.8) to the case of associative algebras. The semi--infinite $\rm Tor$ functor introduced in Section \ref{stor} is the derived functor of the functor of semiproduct.

Let $M\in \Arm$ be a right graded $A$--module and $M^\prime \in \oppAlm$ a left graded $\oppA$--module. Consider the subspace $M\otimes^NM^\prime$ in the tensor product $M\otimes M^\prime$ defined by
$$
M\otimes^NM^\prime=\{ m\otimes m'\in M\otimes M^\prime:~mn\otimes m'=m\otimes nm'\mbox{ for every }n\in N\}.
$$

\begin{definition}\label{spr}
The semiproduct $M\spr M^\prime$ of modules $M\in \Arm$ and $M^\prime\in \oppAlm$ is the image of the subspace $M\otimes^NM^\prime\subset M\otimes M^\prime$ under the canonical projection $M\otimes M^\prime \ra M\otimes_BM^\prime$,
$$
M\spr M^\prime={\rm Im}(M\otimes^NM^\prime \ra M\otimes_BM^\prime).
$$
\end{definition}

Thus the semiproduct $\spr$ is a mixture of the tensor product $\otimes_B$ over $B$ and of the functor $\otimes^N$ of ``N--invariants''. However the following lemma shows that properties of the semiproduct are rather closely related to those of the usual tensor product (compare with \cite{V}, proof of Theorem 3.7). 

\begin{lemma}\label{propspr}
Let $M\in \Armb$ be a right graded $A$--module, $M^\prime \in \oppAlmb$ a left graded $\oppA$--module and $S_A$ the semiregular bimodule associated to $A$. Then
$$
S_A\spr M^\prime = M^\prime
$$
as a left $\oppA$--module, and 

$$
M\spr S_A = M
$$
as a right $A$--module. 
\end{lemma}
 
\pr
We shall prove that $S_A\spr M^\prime$ is isomorphic to $M^\prime$ as a left $\oppA$--module. The second isomorphism may be established in a similar way.

First we calculate the space $S_A\otimes^NM^\prime$. Using realization (\ref{SA2}) of the semiregular bimodule we have:
$$
S_A\otimes^NM^\prime={\rm hom}_B(A,B)\otimes^NM^\prime.
$$

By Lemma \ref{tenshom} and the definition of the operation $\otimes^N$ we also have the following isomorphism of left $B$--modules:
$$
{\rm hom}_B(A,B)\otimes^NM^\prime={\rm hom}_{\k}(N,B)\otimes^N M^\prime=
{\rm hom}_{N}(N,B\otimes M^\prime),
$$
where $N$ acts on $N$ and $M^\prime$ from the left and the $B$--action is induced by the left regular action of $B$ on itself.

Finally the restriction isomorphism
$$
{\rm hom}_{N}(N,B\otimes M^\prime)=B\otimes M^\prime
$$
yields an isomorphism of left $B$--modules $S_A\otimes^NM^\prime= B\otimes M^\prime$ given by
\begin{equation}\label{step1}
\begin{array}{l}
S_A\otimes^NM^\prime={\rm hom}_B(A,B)\otimes^NM^\prime \ra B\otimes M^\prime , \\
\\
f\otimes m' \mapsto f(1)\otimes m'.
\end{array}
\end{equation}

Next we describe the image of the space $S_A\otimes^NM^\prime$ under the canonical projection $S_A\otimes^NM^\prime \ra S_A\otimes_BM^\prime$. Using realization (\ref{SA1}) of the semiregular bimodule and the triangular decomposition $A=N\otimes B$ for the algebra $A$ (see Section \ref{setup}) we obtain the following isomorphism of left $N$--modules
\begin{equation}\label{iso*}
\begin{array}{l}
N^*\otimes M^\prime \ra  N^*\otimes_NA\otimes_B M^\prime=S_A\otimes_B M^\prime, \\
\\
f\otimes m'\mapsto f\otimes 1 \otimes m',
\end{array}
\end{equation}
where the $N$--module structure is induced by the natural left action of $N$ on $N^*$.

Now observe that the isomorphism $\phi:N^*\otimes_N A \ra {\rm hom}_B(A,B)$ constructed in Lemma \ref{phiiso} sends elements of the form $f\otimes 1 \in N^*\otimes_N A,~f\in N^*$ to homomorphisms which take scalar values when restricted to the subspace $N\subset A$. Therefore, recalling isomorphism (\ref{step1}), one can establish an isomorphism of left $B$--modules,
\begin{equation}\label{twostar}
S_A\spr M^\prime = M^\prime.
\end{equation}

Using again isomorphism (\ref{step1}) and observing that (\ref{iso*}) is  an isomorphism of left $N$--modules we conclude that (\ref{twostar}) is in fact an isomorphism of left $N$--modules as well. This completes the proof.

\qed

In conclusion we remark that the semiproduct of modules naturally extends to a functor $\spr:(\Arm ) \times (\oppAlm ) \ra {\rm Vect}_\k$.


\subsection{Semi--infinite homological algebra}\label{shomalg}

\setcounter{equation}{0}
\setcounter{theorem}{0}

In this section we recall, following \cite{V}, the main theorem of semi--infinite homological algebra. Using this theorem we define the semi--infinite Tor functor in the next section. 

The main theorem of semi--infinite homological algebra asserts that the category $\DAlb$ introduced in Section \ref{setup} is equivalent to the homotopy category of semijective complexes.  We remark that in \cite{V} this equivalence was established in the Lie algebra case, $A=\Ug$, where $\g$ is a $\mathbb Z$--graded Lie algebra with finite--dimensional graded components. But in fact the formulation and the proof of the main theorem of semi--infinite homological algebra only use general homological constructions and properties (i),(ii),(v) and (vi) of algebra $A$ axiomatized in Section \ref{setup}. Therefore in this section we reformulate results of \cite{V} for the
 algebra $A$ without any additional comments. 

First we recall (see \cite{V}, Definition 3.3) that a complex $S^{\gr}\in \CAlb$ is called semijective if
{\em 
\vskip 0.3cm
\qquad (1) $S^{\gr}$ is K-injective as a complex of $N$--modules (see \cite{Sp}, Sect. 1), i.e., for every acyclic complex $A^{\gr}\in {\rm Kom}(N)_0$, ${\rm Hom}_{K(N)_0}(A^{\gr},S^{\gr})=0$;
\vskip 0.3cm
\qquad (2) $S^{\gr}$ is K--projective relative to $N$, i.e., for every complex $A^{\gr}\in \CAlb$, such that $A^{\gr}$ is isomorphic to zero in the category $K(N)_0$, ${\rm Hom}_{K(A)_0}(S^{\gr},A^{\gr})=0$.
\vskip 0.3cm
}
By (\ref{homK1}) these two conditions are equivalent to the following ones:
{\em 
\vskip 0.3cm
\qquad (1) For every acyclic complex $A^{\gr}\in {\rm Kom}(N)_0$, the complex  ${\rm hom}_{N}^{\gr}(A^{\gr},S^{\gr})=$ is acyclic;
\vskip 0.3cm
\qquad (2) For every complex $A^{\gr}\in \CAlb$, such that $A^{\gr}$ is homotopic to zero as a complex of $N$--modules, the complex  ${\rm hom}_{A}^{\gr}(S^{\gr},A^{\gr})$ is acyclic.
\vskip 0.3cm
}

Similar definitions may be given for complexes of right $A$--modules from the category $\Armb$.

The main difficulty in dealing with semijective complexes is that in general position the complex of semijective modules is not semijective. However in some particular cases described in the next proposition K--injectivity (K--projectivity relative to $N$ or semijectivity) of the complex follows from the corresponding property of the individual terms of this complex.
\begin{proposition}{\bf (\cite{V}, Proposition 3.7)}\label{sinjprop}

\noindent
1. Any complex $S^{\gr}\in \CAlb$ of $N$--injective modules bounded from below is K--injective as a complex of $N$--modules.

\noindent
2. Any complex $S^{\gr}\in \CAlb$ of projective relative to $N$ modules bounded from above is K--projective relative to $N$.  

\noindent
3. Any bounded complex $S^{\gr}\in \CAlb$ of semijective modules is semijective.
\end{proposition}

An $A$--module $M\in \Almb$ is called semijective if the corresponding 0--complex $\ldots\ra 0\ra M\ra 0\ra \ldots$ (see Section \ref{notat}) is semijective.
We also say that $M$ is projective relative to $N$ if the corresponding 0--complex is K--projective relative to $N$. For the 0--complex $\ldots\ra 0\ra M\ra 0\ra \ldots$ condition 1 of the definition of semijective complexes is equivalent to the usual $N$-injectivity of $M$.

In this paper we shall actually deal with a class of relatively to $N$ projective modules described in the next lemma (see \cite{V}, Sect. 3.1).
\begin{lemma}\label{relproj}
Every left $A$--module $M\in \Almb$ induced from an $N$--module $V\in (N-{\rm mod})_0$, $M=A\otimes_NV$, is projective relative to $N$.
\end{lemma}
\pr
Let $A^{\gr}$ be a complex of left $A$--modules from the category $\Armb$ such that $A^{\gr}$ is homotopic to zero as a complex of $N$--modules. We have to show that $H^\gr({\rm hom}_{A}^{\gr}(M,A^\gr))=0$.

Indeed, since $M=A\otimes_NV$ we have a canonical isomorphism of complexes ${\rm hom}_{A}^{\gr}(M,A^{\gr})={\rm hom}_{N}^{\gr}(V,A^{\gr})$. But the complex ${\rm hom}_{N}^{\gr}(V,A^{\gr})$ is homotopic to zero and, in particular, acyclic since $A^{\gr}$ is homotopic to zero as a complex of $N$--modules. This completes the proof.

\qed

The following fundamental property of the semiregular bimodule $S_A$ together with Lemma \ref{propspr} shows that $S_A$ is an analogue of the regular representation in semi--infinite homological algebra. 
\begin{proposition}\label{SAprop}
 Let $A$ be  an associative $\mathbb Z$--graded algebra over a field $\k$ satisfying conditions (i)--(viii) of Sections \ref{setup} and  \ref{bimod}.  
Then the semiregular bimodule $S_A$ is semijective as a right $A$--module and a left $\oppA$--module.
\end{proposition}
\pr(Compare with \cite{V}, Proposition 3.6)
Consider $S_A$ as a right $A$--module. 
First we prove that $S_A$ is injective as a right $N$--module.
Using realization (\ref{SA2}) of the right semiregular representation and property (vi) of the algebra $A$ we obtain the following isomorphism of right $N$--modules $S_A={\rm hom}_B(A,B)={\rm hom}_{\k}(N,B)$. The last module is evidently $N$--injective.

It is also clear that $S_A$ is projective relative to $N$. Indeed,  from realization (\ref{SA1}) of the right semiregular representation we obtain that as a right $A$--module $S_A$ is induced from $N$--module $N^*$, $S_A=N^*\otimes_NA$, and hence by Lemma \ref{relproj} $S_A$ is projective relative to $N$ as a right $A$--module. We conclude that $S_A$ is semijective as a right $A$--module.

The proof of the fact that $S_A$ is semijective as a left $\oppA$--module is quite similar to the one presented above. We just have to use two realizations of the left semiregular action of $\oppA$ on $S_A$ obtained in Proposition \ref{SAopp}.

\qed 

Now we formulate the main theorem of semi--infinite homological algebra.
\begin{theorem}{\bf (\cite{V}, Theorem 3.3)}\label{mainsinf}
Let $A$ be an associative $\mathbb Z$--graded algebra satisfying conditions (i),(ii),(v) and (vi) of Section \ref{setup}. Let ${\rm Kom}({\mathcal S \mathcal J}(A)_0)$ be the category of semijective complexes associated to the abelian category $\Almb$.
Denote by  $K({\mathcal S \mathcal J}(A)_0)$ the corresponding homotopy category. Then the functor of localization by the class of quasi--isomorphisms is an equivalence of categories:
$$
K({\mathcal S \mathcal J}(A)_0)\cong \DAlb.
$$
\end{theorem}

In particular, we have the following important corollary of Theorem \ref{mainsinf}.

\begin{corollary}{\bf (\cite{V}, Theorem 3.2)}\label{sires}
For every complex $K^{\gr}\in \CAlb$ there exists a quasi--isomorphism $S^{\gr}\ra K^{\gr}$, where $S^{\gr}\in \CAlb$ is a semijective complex. The complex $S^{\gr}$ is called a semijective resolution of $K^{\gr}$. 
\end{corollary}

Properties of semijective resolutions are summarized in the following proposition that is also a corollary of Theorem \ref{mainsinf}.

\begin{proposition}{\bf (\cite{V}, Corollaries 3.1 and  3.2)}\label{siresprop}
Let $\phi:K^{\gr}\ra {K^\prime}^{\gr}$ be a morphism in $\DAlb$, and $S^{\gr},~{S^\prime}^{\gr}$ semijective resolutions of $ K^{\gr}$ and ${K^\prime}^{\gr}$, respectively. Then there exists a morphism of complexes $\phi^{\gr}: S^{\gr}\ra {S^\prime}^{\gr}$ in the category $\CAlb$ such that the square
$$
\begin{array}{ccc}
S^{\gr}~~~~~~ & \longrightarrow & K^{\gr} \\
\downarrow {\phi^{\gr}}& ~~~ & \downarrow {\phi} \\
{S^{\prime}}^{\gr}~~~~~~ & \longrightarrow & {K^\prime}^{\gr}
\end{array}
$$
is commutative in $\DAlb$. This morphism is unique up to a homotopy. 

In particular, any two semijective resolutions of a complex $K^{\gr}$ are homotopically equivalent. This equivalence is unique up to a homotopy.
\end{proposition}

\begin{corollary}{\bf (\cite{V}, Corollary 3.3)}\label{acycl}
Each acyclic semijective complex is homotopic to zero. 
\end{corollary}

By definition a semijective resolution of a left $A$--module $M\in \Almb$ is a semijective resolution of the corresponding 0--complex $\ldots\ra 0 \ra M\ra 0 \ra \ldots$. Next we formulate, for future references, properties of semijective resolutions of left $A$--modules. These properties follow directly from Corollary \ref{sires} and Proposition \ref{siresprop} applied to 0--complexes $K^\gr= \ldots\ra 0 \ra M\ra 0 \ra \ldots$.
\begin{proposition}\label{sres}

(a) Every left $A$--module $M\in \Almb$ has a semijective resolution.

(b) Any morphism of $A$--modules $M,~M^\prime \in \Almb$, $\phi:M\ra M^\prime$, gives rise to a morphism (in the category $\CAlb$) of their semijective resolutions $\phi^{\gr}: S^{\gr}\ra {S^{\prime}}^\gr$, that is unique up to a homotopy. 

(c) In particular, any two semijective resolutions of a module $M\in \Almb$ are homotopically equivalent. This equivalence is unique up to a homotopy.
\end{proposition}

Using Proposition \ref{ocompl} we obtain the following simple characterization of semijective resolutions of modules (see \cite{V}, Sect. 3.4).
\begin{proposition}\label{resmod}
A semijective complex $S^{\gr}\in \CAlb$ is a semijective resolution of a module $M\in \Almb$ if and only if
$$
H^i(S^{\gr})=\left\{
\begin{array}{l}
M~\mbox{ for }i=0 \\
0~\mbox{ for }i\neq 0
\end{array}
\right
 .
$$
\end{proposition}

In conclusion we note that the results of this section may be carried over to the homotopy and derived categories associated to the abelian categories $\Armb$, $\oppAlmb$ and $\oppArmb$.


\subsection{Semi--infinite Tor functor}\label{stor}

\setcounter{equation}{0}
\setcounter{theorem}{0}

In this section we define, using the results of the previous section, the semi--infinite Tor functor as the classical derived functor of the functor of semiproduct introduced in Section \ref{sspr}. We suppose that the algebra $A$ satisfies conditions (i)--(viii) of Sections \ref{setup} and \ref{bimod}.

First we define the semi--infinite Tor functor on modules $M\in \Armb$, $M^\prime \in \oppAlmb$ as the cohomology space of the complex $S^{\gr}(M)\spr S^{\gr}(M^\prime)$,
$$
\stor(M,M^\prime)=H^{\gr}(S^{\gr}(M)\spr S^{\gr}(M^\prime)),
$$ 
where $S^{\gr}(M),~S^{\gr}(M^\prime)$ are semijective resolutions of $M$ and $M^\prime$. By Proposition \ref{sres} (c) the space $\stor(M,M^\prime)\in {\rm Kom}({\rm Vect}_{\k})$ does not not depend on the resolutions $S^{\gr}(M),~S^{\gr}(M^\prime)$. We shall show that $\stor$ may be naturally extended to a functor $\stor: \Armb \times \oppAlmb \ra {\rm Kom}({\rm Vect}_{\k})$.

Indeed, consider two morphisms of modules $\phi: M\ra \tilde M$ and $\phi^\prime: M^\prime \ra \tilde M^\prime$, where $M,~\tilde M \in \Armb$ and  $M^\prime,~\tilde M^\prime \in \oppAlmb$. Using  part (b) of Proposition \ref{sres} these morphisms of modules  give rise to  morphisms $\phi^{\gr}: S^{\gr}(M) \ra S^{\gr}(\tilde M)$ and ${\phi^\prime}^{\gr}: S^{\gr}(M^\prime) \ra S^{\gr}(\tilde M^\prime)$ of semijective resolutions $S^{\gr}(M)$, $S^{\gr}(\tilde M)$, $S^{\gr}(M^\prime)$ and $S^{\gr}(\tilde M^\prime)$ of these modules which are unique up to homotopies. Therefore  one can define a natural map
$$
\stor(\phi,\phi^\prime): \stor(M,M^\prime)\ra \stor(\tilde M,\tilde M^\prime).
$$
We conclude that modules $\stor(M,M^\prime)$ together with maps $\stor(\phi,\phi^\prime)$ yield a functor 
$$
\stor: \Armb \times \oppAlmb \ra {\rm Kom}({\rm Vect}_{\k}).
$$
This functor is called the semi--infinite Tor functor.

The following important theorem is a semi--infinite analogue of the classical theorem about partial derived functors (see \cite{carteil}, Ch. V, \S 8, Theorem 8.1).
\begin{theorem}\label{threetor}
The following three definitions of the spaces \\ $\stor(M,M^\prime)\in {\rm Kom}({\rm Vect}_{\k})$ are equivalent:
\vskip 0.3cm
(a) $\stor (M,M^\prime)=H^{\gr}(S^{\gr}(M)\spr S^{\gr}(M^\prime))$;
\vskip 0.3cm
(b) $\stor (M,M^\prime)=H^{\gr}(M\spr S^{\gr}(M^\prime))$;
\vskip 0.3cm
(c) $\stor (M,M^\prime)=H^{\gr}(S^{\gr}(M)\spr M^\prime)$,
\vskip 0.3cm
\noindent
where $M\in \Armb$, $M^\prime \in \oppAlmb$, and $S^{\gr}(M),~S^{\gr}(M^\prime)$ are semijective resolutions of $M$ and $M^\prime$, respectively.
\end{theorem} 

\begin{remark}
Clearly, beside of the functor $\stor$ one can introduce partial derived functors of the functor of semiproduct. On objects 
$M\in \Armb$, $M^\prime \in \oppAlmb$ they are given by formulas (b) and (c) of the previous theorem, on morphisms they are defined similarly to the maps $\stor(\phi,\phi^\prime)$ (see the definition of the functor $\stor$).
Theorem \ref{threetor} only establishes isomorphisms of the complexes of graded vector spaces defined by formulas (a), (b) and (c). In contrast to the classical case we  have no isomorphisms of the corresponding derived functors. 
\end{remark} 

If one of two modules $M\in \Armb$, $M^\prime \in \oppAlmb$ is semijective, then this module is a semijective resolution of itself. Therefore we have the following simple corollary of Theorem \ref{threetor}.
\begin{corollary}\label{svanishtor}
Suppose that one of modules $M\in \Armb$, $M^\prime \in \oppAlmb$ is semijective. Then
$$
\stor (M,M^\prime)=M\spr M^\prime.
$$
\end{corollary}

The proof of Theorem \ref{threetor} occupies Appendix 5. In this proof we  use two types of standard semijective resolutions constructed in the next section.


\subsection{Standard semijective resolutions}\label{stres}

\setcounter{equation}{0}
\setcounter{theorem}{0}

In this section, using standard relative bar resolutions, we construct two types of standard semijective resolutions. We start by recalling the definition of the standard (normalized) relative bar resolution (see \cite{guish}, Appendix C and \cite{Arkh1}, Sect. 2.2).

Let $A$ be a $\mathbb Z$--graded associative algebra over a field $\k$ satisfying conditions (i) and (vi) of Section \ref{setup}. The standard bar resolution $\tilBar (A,B,M)$ of a left $\mathbb Z$--graded $A$-module $M$ with respect to the subalgebra $B\subset A$ is defined as follows:
\begin{equation}\label{bar}
\begin{array}{l}
  \tilBarn (A,B,M)=\underbrace{A\otimes_B\ldots \otimes_BA}_{n+1 ~\mbox{\tiny  times}}\otimes_BM,~~n\leq 0,  \\
 \\
 d(a_0\otimes \ldots \otimes a_n\otimes v)= \\
\\
  \sum_{s=0}^{n-1}(-1)^s a_0\otimes\ldots\otimes
  a_sa_{s+1}\otimes\ldots\otimes v + \\
 \\
  +(-1)^na_0\otimes\ldots\otimes a_{n-1}\otimes a_nv,
\end{array}
\end{equation}
where
$a_0,\ldots ,a_n\in A,\ v\in M$.

In order to define the standard normalized relative bar resolution we need the following simple lemma.
\begin{lemma}{\bf (\cite{Arkh1}, Lemma 2.2.1)}
The subspace $\linBar (A,B,M)$,
$$
\begin{array}{l}
\linBarn(A,B,M)= \\
 
\{ a_0\otimes\ldots\otimes a_n\otimes
  v\in \tilBarn(A,B,M)|~\exists s\in\{1,\ldots,n\}:
  a_s\in B\}
\end{array}
$$
is a subcomplex in $\tilBar (A,B,M)$.
\end{lemma}

The quotient complex $\Bar (A,B,M)=\tilBar (A,B,M)/\linBar
(A,B,M)$ is called the normalized bar resolution of the
$A$--module $M$ with respect to the subalgebra $B$.

The following  properties of the standard normalized bar resolution may be checked directly using definition (\ref{bar}) of the complex $\tilBar(A,B,M)$.

\begin{proposition}\label{barprop}
Let $M$ be a left $\mathbb Z$--graded $A$--module. Then
\vskip 0.3cm
(i)$\Bar (A,B,M)$ is a resolution of $M$, i.e. the natural map $\Bar (A,B,M)\ra M$ is a quasi--isomorphism.
\vskip 0.3cm
(ii) $\Bar (A,B,M)= \Bar (N,{\k} ,M)$ as a complex of $N$--modules. In particular $\Bar (A,B,M)$ is an $N$--free resolution of $M$.
\vskip 0.3cm
(iii) The complex $\Bar (A,B,A)$ is homotopically equivalent to $A$ as a complex of $A$--$B$ and $B$--$A$ bimodules, the bimodule structures being induced by the left and right regular actions of $A$. The corresponding homotopy maps are given by
$$
a_0\otimes\ldots\otimes a_n\otimes a \mapsto a_0\otimes\ldots\otimes a_n\otimes a \otimes 1
$$
and
$$
a_0\otimes\ldots\otimes a_n\otimes a \mapsto 1\otimes a_0\otimes\ldots\otimes a_n\otimes a,
$$
respectively.
\end{proposition}

Now using normalized bar resolutions we construct two types of standard semijective resolutions of modules.

\begin{proposition}\label{res1}
Let $M\in \Armb$ be a right $A$-module. Then the complex $\sBar(A,N,M)$  defined by
$$
\sBar(A,N,M)=\hAd(\Bar(A,B,A),M)\otimes_A\Bar(A,N,A)
$$
is a semijective resolution of $M$.
\end{proposition}

\begin{proposition}\label{res2}
Let $M^\prime \in \oppAlmb$ be a left $\oppA$--module. Then the complex $\sBaropp$ defined by
$$
\sBaropp(\oppA,N,M^\prime)=\sBar(A,N,S_A)\spr M^\prime
$$
is a semijective resolution of $M^\prime$.
\end{proposition}

The proofs of these two propositions are contained in Appendixes 3 and 4, respectively.
Here we only prove the following weak version of Proposition \ref{res1}.
\begin{proposition}\label{weak}
Let $M\in \Armb$ be a right $A$-module. Suppose that $M$ is injective as an $N$--module. Then the complex 
\begin{equation}\label{stresinj}
M\otimes_A\Bar(A,N,A)
\end{equation}
is a semijective resolution of $M$.
\end{proposition}
\pr
By Proposition \ref{resmod} we have to prove that $M\otimes_A\Bar(A,N,A)$ is a semijective complex such that $H^\gr(M\otimes_A\Bar(A,N,A))=M$.

First observe that by the definition of the normalized bar resolution the natural map $M\otimes_A\Bar(A,N,A)\ra M$ is quasi--isomorphism. Therefore $H^\gr(M\otimes_A\Bar(A,N,A))=M$.

Next we show that the complex $M\otimes_A\Bar(A,N,A)$ is $N$--K--injective.
Indeed, by part (iii) of Proposition \ref{barprop} the complex $\Bar(A,N,A)$ is homotopically equivalent to $A$ as a complex of $A$--$N$--bimodules, and hence we have a homotopy equivalence of complexes of $N$--modules,
$$
M\otimes_A\Bar(A,N,A)\cong \ldots \ra 0\ra M\ra 0 \ra \ldots.
$$
The complex $\ldots \ra 0\ra M\ra 0 \ra \ldots$ is $N$--K--injective as a bounded complex of $N$--injective modules (see part 3 of Proposition \ref{sinjprop}). We conclude that the complex $M\otimes_A\Bar(A,N,A)$ is also $N$--K--injective.

Now we prove that the complex $M\otimes_A\Bar(A,N,A)$ is relatively to $N$ K--projective.
From the definition of the standard normalized bar resolution it follows that the individual terms of the complex $M\otimes_A\Bar(A,N,A)$ are $A$--modules induced from $N$--modules. Therefore by Lemma \ref{relproj} the individual terms of the complex $M\otimes_A\Bar(A,N,A)$ are relatively to $N$ projective modules.

Now by part 2 of Proposition \ref{sinjprop} the  complex $M\otimes_A\Bar(A,N,A)$ is  relatively to $N$ K--projective as a bounded from above complex of relatively to $N$ projective modules. This completes the proof.

\qed

A resolution similar to (\ref{stresinj}) may be constructed for left $N$--injective $\oppA$--modules. 

Substituting the resolution obtained in Proposition \ref{res1} and the modification of this resolution for left $\oppA$--modules into formulas (b) and (c) of Theorem \ref{threetor}, respectively, we obtain the following corollary of Proposition \ref{weak}.
\begin{corollary}
Suppose that one of modules $M\in \Armb$, $M^\prime \in \oppAlmb$ is $N$--injective. Then
$$
{\rm Tor}_A^{\frac{\infty}{2}+0}(M,M^\prime)=M\spr M^\prime,
$$
and
$$
{\rm Tor}_A^{\frac{\infty}{2}+n}(M,M^\prime)=0~ \mbox{ for } n>0.
$$
\end{corollary}


\subsection{Semi--infinite cohomology of Lie algebras}\label{lie}

\setcounter{equation}{0}
\setcounter{theorem}{0}

In this section we discuss the semi--infinite Tor functor for the class of $\mathbb Z$--graded Lie algebras with finite--dimensional graded components. 

Let $\g$ be a $\mathbb Z$--graded Lie algebra over $\k$, $\g=\bigoplus_{n\in {\mathbb Z}}\g_n$, such that ${\rm dim}~\g_n<\infty$. Denote the subalgebras $\bigoplus_{n>0}\g_n$ and $\bigoplus_{n\leq 0}\g_n$ by $\n$ and $\b$, respectively. Let $\Ug,~\Ub$ and $\Un$ be the universal enveloping algebras of $\g$ and of the subalgebras $\b,\n\subset \g$. The  algebras $A=\Ug,~B=\Ub$ and $N=\Un$ satisfy conditions (i)--(viii) of Sections \ref{setup}, \ref{bimod} (see \cite{Arkh1}, Sect. 4), and hence one can define the algebra $\oppUg$ and  the semi--infinite Tor functor for $\Ug$.
Remarkably, the algebra $\oppUg$ may be described explicitly as the universal enveloping algebra of the central extension of $\g$ with the help of the critical two--cocycle. 

The critical cocycle on $\g$ may be defined as follows. (see \cite{V}, Sect. 2).
First consider a $\mathbb Z$--graded vector space $V$ over a field $\k$, $V=\bigoplus_{n\in {\mathbb Z}}V_n$ such that ${\rm dim}~V_n<\infty$ for $n>0$. Denote the spaces $\bigoplus_{n>0}V_n$ and $\bigoplus_{n\leq 0}V_n$ by $V_+$ and $V_-$, respectively. 
We shall need a certain subalgebra, that we denote by ${\g\mathfrak l}(V)$, in the algebra of linear transformations of the space $V$. To define this subalgebra
we observe that every linear transformation $a$ of the space $V$ may be represented in the form 
$$
a=\left(\begin{array}{cc}a_{++}&a_{++}\\ a_{-+}&a_{--}\end{array}\right),
$$
where $a_{++}:V_+\ra V_+$, $a_{+-}:V_-\ra V_+$, etc., are the blocks of $a$ with respect to the decomposition $V=V_+\oplus V_-$ of the space $V$.

Now we define ${\g\mathfrak l}(V)$ as the subalgebra in the algebra  of linear transformations of $V$ formed by elements $a$ such that their blocks $a_{-+}:V_+\ra V_-$ are of finite rank. 

The algebra ${\g\mathfrak l}(V)$ has a remarkable central extension with the help of two--cocycle $\omega_V$ (see, for instance, \cite{V}, Proposition 2.1),
$$
\omega_V(a,b)={\rm tr}(b_{-+}a_{+-}-a_{-+}b_{+-}).
$$

Applying the construction of the cocycle $\omega_V$ for $V=\g$ we obtain a cocycle $\omega_\g$ defined on the Lie algebra ${\g\mathfrak l}(\g)$.  
The adjoint representation of $\g$ provides a morphism of Lie algebras $\g \ra  {\g\mathfrak l}(\g)$. The inverse image of the two--cocycle $\omega_\g$ under this morphism is called the critical cocycle of $\g$. We denote this cocycle by $\omega_0$.

Now we can explicitly describe the algebra $\oppUg$. 
\begin{proposition}{\bf (\cite{Arkh1}, Corollary 4.4.2)}\label{oppug}
Let $\oppg=\g+K\k$ be the central extension of $\g$ with the help of the critical cocycle $-\omega_0$. Then the algebra $\oppUg$ is isomorphic to the quotient $U(\oppg)/I$, where $I$ is the two--sided ideal in $U(\oppg)$ generated by $K-1$.
\end{proposition}

The semi--infinite cohomology was first defined for the class of graded Lie algebras equipped with so--called semi--infinite structures (see \cite{F,V}). Below we show that the semi--infinite cohomology spaces proposed in \cite{F,V} may be defined using the semi--infinite Tor functor for associative algebras (see Section \ref{stor}). 

First we recall the definition of the semi--infinite structure for the Lie algebra $\g$.
Suppose that the cohomology class of the critical cocycle $\omega_0$ is trivial in $H^2(\g,\k)$. One says that the Lie algebra $\g$ is equipped with a semi--infinite structure if there exists a one--cochain $\beta\in \g^*$ on $\g$ such that $\beta=0$ on $\g_n$ for all $n\neq 0$, and $\partial \beta=\omega_0$.
        
Observe that the cochain $\beta:\g \ra \k$ naturally extends to a one--dimensional representation of the algebra $\oppg$, $\beta(x,k)=\beta(x)+k$, where $(x,k)\in \g+K\k =\oppg$. Indeed, it suffices to verify that  $\beta([(x,k),(x',k')])=0$ for any $(x,k),(x',k')\in \oppg$. Using Proposition \ref{oppug} and the definition the cochain $\beta$ we have
$$
\begin{array}{l}
\beta([(x,k),(x',k')]) = \beta([x,x'],-\omega_0(x,x'))= \\
\\
\beta([x,x'])-\omega_0(x,x')=\omega_0(x,x')-\omega_0(x,x')=0.
\end{array}
$$

The representation $\beta:\oppg \ra \k$ naturally extends to a one--dimensional representation of the universal enveloping algebra $U(\oppg)$.
Since the two--sided ideal in $U(\oppg)$ generated by $K-1$ lies in the kernel of the representation $\beta: U(\oppg)\ra \k$ this representation gives rise to
a representation of the algebra $\oppUg$.
We denote this one--dimensional left $\oppUg$--module by $\kb$. 

Let $M\in \Ugrmb$ be a right $\Ug$--module. The semi--infinite cohomology space of $M$ is defined as
$$
H^{\frac{\infty}{2}+\gr}(\Ug,M)=\storg(M,\kb).
$$

The properties of the simi--infinite cohomology spaces $H^{\frac{\infty}{2}+\gr}(\Ug,M),$ $M\in \Ugrmb$ may be derived from the general properties of the semi--infinite Tor functor (see Section \ref{stor}). 
In particular, part (c) of Theorem \ref{threetor} shows that the spaces $H^{\frac{\infty}{2}+\gr}(\Ug,M)$ coincide with the semi--infinite cohomology spaces defined in \cite{V}, Sect. 3.9 as 
$$
H^{\frac{\infty}{2}+\gr}(\Ug,M)=H^\gr(S^\gr(M)\otimes_{\Ub}^{\Un} \kb),
$$
where $S^\gr(M)$ is a semijective resolution of $M$.

However in \cite{V} Voronov  also proves another interesting vanishing theorem for the spaces $H^{\frac{\infty}{2}+\gr}(\Ug,M)$ that we still could not verify in case of arbitrary associative algebras.

\begin{proposition}{\bf (\cite{V}, Theorem 2.1)}\label{vv}
Let $M\in \Ugrmb$ be a right $\Ug$--module. Suppose that $M$ is injective as an $\Un$--module and projective as an $\Ub$--module. Then
$$
H^{\frac{\infty}{2}+\gr}(\Ug,M)=M\otimes_{\Ub}^{\Un} \kb.
$$
\end{proposition}

Next following \cite{V} we shall describe the standard Feigin's complex for calculation of the semi--infinite cohomology of $M$.
First we recall the construction of the standard semijective resolution of the one--dimensional $\oppUg$--module $\kb$ (see \cite{V}, Sect. 3.7). This resolution is a semi--infinite analogue of the standard resolution $\Ug\otimes \Lambda^{\gr}(\g)$  of the trivial $\Ug$--module (see \cite{carteil}, Ch.XIII, \S 7), the regular bimodule $\Ug$ and the exterior algebra $\Lambda^{\gr}(\g)$ being replaced with their semi--infinite counterparts $S_{\Ug}$ and $\sl$, where $\sl$ is the space of semi--infinite exterior forms  on $\g^*$. This space is defined as follows.
 
Consider the Clifford algebra $C(\g+\g^*)$ associated  with the symmetric bilinear form
$$
\begin{array}{l}
\langle v_1^*,v_2\rangle =v_1^*(v_2), \\
\\
\langle v_1,v_2\rangle =\langle v_1^*,v_2^*\rangle =0
\end{array}
$$
on the vector space $\g +\g^*$, where $\g^*=\bigoplus_{n\in {\mathbb Z}}\g_n^*$,
$v_{1,2}\in \g,~~v_{1,2}^*\in \g^*$. For every $v\in \g,~v^*\in \g^*$ we denote by $\overline v,~\overline v^*$ the elements of the algebra $C(\g +\g^*)$ which correspond to $v$ and $v^*$, respectively.
The algebra $C(\g+\g^*)$ is generated by the vector spaces $\g$ and $\g^*$, with the defining relations
\begin{equation}\label{cl}
\overline a\overline b+\overline b\overline a=\langle a,b\rangle , \mbox{ where } a,b \in \g+\g^*.
\end{equation}

The space $\sl$  of semi--infinite exterior forms on $\g^*$ is defined as the representation of the algebra $C(\g+\g^*)$ freely generated by the vacuum vector $x_0$ that satisfies the conditions
$$
\begin{array}{l}
\overline v\cdot x_0=0 \mbox{ for } v\in \n, \\
\\
\overline v^*\cdot x_0=0 \mbox{ for } v^*\in \n^\perp = \b^*.
\end{array}
$$

Choose a linear basis $\{ e_i\}_{i\in {\mathbb Z}}$ of $\g$ compatible with the $\mathbb Z$ grading on $\g$  in the sense that for any $i\in {\mathbb Z}$ $e_i\in \g_n$ for some $n\in {\mathbb Z}$ and if  $e_i\in \g_n$ then $e_{i+1}\in \g_n$ or $e_{i+1}\in \g_{n+1}$. Let $\{ e_i^*\}_{i\in {\mathbb Z}}$ be the dual basis of $\g^*$.
Each element of $\sl$ is a linear combination of monomials of the type
\begin{equation}\label{mon}
\omega = \overline e_{i_1}\ldots \overline e_{i_m}\overline e_{j_1}^*\ldots \overline e_{j_n}^*\cdot x_0.
\end{equation}
From commutation relations (\ref{cl}) it follows that for each monomial $\omega$ of the form (\ref{mon}) the integer number ${\rm deg}~\omega = n-m$ is well defined. This number is called the degree of $\omega$. This equips the space $\sl$ with a $\mathbb Z$--grading. 

Remark that if $\g=\b$ the space $\sl$ degenerates into the exterior algebra $\Lambda^\gr(\g)$.

Now we describe the standard semijective resolution of the one--dimensional $\oppUg$--module $\kb$.
Equip the left $\oppUg$--module $S_{\Ug}\otimes \sl$ with an operator $d$,
\begin{equation}\label{different}
d=\sum_i e_i\otimes \overline{e}_i^*-\sum_{i<j} 1\otimes :\overline{[e_i,e_j]} \overline e_i^* \overline e_j^*: +1\otimes \overline \beta.
\end{equation}
Here $:\overline{[e_i,e_j]} \overline e_i^* \overline e_j^*:$ is the normally ordered product of the elements $\overline{[e_i,e_j]}, \overline e_i^*, \overline e_j^*$, i.e. $:\overline{[e_i,e_j]} \overline e_i^* \overline e_j^*:$ is a permuted product of the elements $\overline{[e_i,e_j]}, \overline e_i^*, \overline e_j^*$ in the algebra $C(\g+\g^*)$, such that all the operators annihilating the vacuum vector $x_0$ of the representation $\sl$ stand on the right, times the sign of the permutation; $e_i\otimes 1$ is regarded as the 
operator of right semiregular action of $\Ug$ on $S_{\Ug}$, $\overline\beta,\overline e_i^*,\overline e_j^*\in C(\g+\g^*)$ act in the space $\sl$ of semi--infinite exterior forms on $\g^*$. 

Note that the normal product operation is well defined because the subspace of operators in $C(\g+\g^*)$ annihilating the vacuum vector $x_0$ is generated by elements from the maximal isotropic subspace  $\n+\b^*\subset \g+\g^*$. 

The operator $d$ has degree 1 with respect to the $\mathbb Z$--grading of $S_{\Ug}\otimes \sl$ by degrees of semi--infinite exterior forms.

The space $S_{\Ug}\otimes \sl$ inherits also the second $\mathbb Z$--grading from the $\mathbb Z$--grading of the Lie algebra $\g$, 
$$
\begin{array}{l}
{\rm deg}'(u\otimes \overline e_{i_1}\ldots \overline e_{i_m}\overline e_{j_1}^*\ldots \overline e_{j_n}^*\cdot x_0)= \\
\\
{\rm deg}'(u)+\sum_{k=1}^m{\rm deg}'(e_{i_k}) - \sum_{l=1}^m{\rm deg}'(e_{j_l}),
\end{array}
$$
where $u\in S_{\Ug}$, $\overline e_{i_1}\ldots \overline e_{i_m}\overline e_{j_1}^*\ldots \overline e_{j_n}^*x_0\in \sl$, ${\rm deg}'(u)$ is the degree of $u$ with respect to the natural $\mathbb Z$--grading in $S_{\Ug}$, and ${\rm deg}'(e_{i_k}), {\rm deg}'(e_{j_l})$ are the degrees of the elements $e_{i_k},~e_{j_l}$ in $\g$. 

Since $\beta \in \g_0^*$ the operator $d$ preserves the second $\mathbb Z$--grading of $S_{\Ug}\otimes \sl$. Moreover, we have the following proposition.
\begin{proposition}{\bf (\cite{V}, Propositions 2.6, 3.13)}\label{stresU}

(i) $d^2=0$, i.e. $d$ equips the space $S_{\Ug}\otimes \sl$ with the structure of a complex.
 
(ii)  The complex $S_{\Ug}\otimes \sl\in \CoppUglb$ is a semijective resolution of the one--dimensional $\oppUg$--module $\kb$.
\end{proposition}

Using part (b) of Theorem \ref{threetor} we obtain the following corollary of Proposition \ref{stresU}.
\begin{corollary}\label{feigin}
Let $M\in \Ugrmb$ be a right $\Ug$--module. The semi--infinite cohomology space of $M$,
$$
H^{\frac{\infty}{2}+\gr}(\Ug,M)=\storg(M,\kb),
$$
may be calculated as the cohomology of the complex $M\otimes_{\Ub}^{\Un} S_{\Ug}\otimes \sl=M\otimes \sl$,
$$
H^{\frac{\infty}{2}+\gr}(\Ug,M)=H^\gr(M\otimes \sl).
$$

The complex $M\otimes \sl$ is called the standard Feigin's complex for calculation of the semi--infinite cohomology of $M$.
\end{corollary}


\section{Semi--infinite Hecke algebras}\label{shecke}

In this section we define semi--infinite modifications of Hecke algebras. Properties of these algebras are quite similar to those of the usual Hecke algebras introduced in Section \ref{hecke}. However, as we shall see in Section \ref{Hklie} the semi--infinite Hecke algebras are rather adapted for the study of ``infinite--dimensional'' objects, e.g. affine Lie algebras.


\subsection{Semi--infinite Hecke algebras: definition and main properties}\label{sheckedefr}

\setcounter{equation}{0}
\setcounter{theorem}{0}

In this section we generalize the notion of Hecke algebras to semi--infinite cohomology. The exposition in this section is parallel to that of Section \ref{hecke}. We also use the notation introduced in Section \ref{setup}.

Let $A$ be an associative $\mathbb Z$ graded algebra over a field $\k$. Suppose that the algebra $A$ contains a graded subalgebra $A_0$, and  both $A$ and $A_0$ satisfy conditions (i)--(viii) of Sections \ref{setup} and \ref{bimod}. We denote by $N,~B$ and $N_0,~B_0$ the graded subalgebras in $A$ and  $A_0$, respectively, providing the triangular decompositions of these algebras (see condition (vi) of Section \ref{setup}). We also assume that the algebra $\oppAo$ is augmented, $\e:\oppAo\ra \k$ (note the difference between this condition and the corresponding condition of Section \ref{hecke}). We denote this one--dimensional $\oppAo$--module by $\ke$.

\begin{definition}
Let $S^\gr({\ke})$ be a semijective resolution of the left $\oppAo$--module $\ke$. Consider the complex $S_A\spro S^\gr({\ke})$ of left $\oppA$--modules, where $A_0$ acts on the semiregular bimodule $S_A$ by right semiregular action and the structure of a left $\oppA$--module on $S_A\spro S^\gr({\ke})$ is induced by the left semiregular action of $\oppA$ on $S_A$. The double graded algebra
\begin{equation}\label{sheckedef}
\Hks=\hDoppAl(S_A\spro S^\gr({\ke}), S_A\spro S^\gr({\ke}))
\end{equation}
is called the semi--infinite Hecke algebra of the triple $(A,A_0,{\e})$.
\end{definition}

From Proposition \ref{ocompl} we obtain the following simple but important vanishing property for semi--infinite Hecke algebras.
\begin{proposition}\label{svanish}
Assume that 
$$
H^\gr(S_A\spro S^\gr({\ke}))=\storo(S_A,\ke)=S_A\spro \ke.
$$ 
Then
$$
\Hks=\hDoppAl(S_A\spro {\ke}, S_A\spro {\ke}).
$$
In particular,
$$
{\rm Hk}^{\frac{\infty}{2}+0}(A,A_0,\e)={\rm hom}_{\oppA}(S_A\spro \ke,S_A\spro \ke).
$$
\end{proposition}

\begin{remark} 
From Corollary \ref{svanishtor} it follows that the condition $\storo(S_A,\ke)=$\\ $S_A\spro \ke$ is satisfied if $S_A$ is semijective as a right $A_0$--module.
\end{remark}

Now similarly to Proposition \ref{hact} we define an action of semi--infinite Hecke algebras in semi--infinite cohomology spaces.  First for every right $A_0$--module $M$, $M\in \Aormb$, we introduce  the semi--infinite cohomology space $H^{\frac{\infty}{2}+\gr}(A_0,M)$ of $M$ by
\begin{equation}\label{shom}
H^{\frac{\infty}{2}+\gr}(A_0,M)=\storo(M,\ke).
\end{equation}

\begin{proposition}
For every right $A$--module  $M\in \Armb$ the algebra \\ $\Hks$ naturally acts in the semi--infinite cohomology space $H^{\frac{\infty}{2}+\gr}(A_0,M)$ of $M$ regarded as a right $A_0$--module,
$$
\Hks \times H^{\frac{\infty}{2}+\gr}(A_0,M)\ra H^{\frac{\infty}{2}+\gr}(A_0,M).
$$
This action respects the bigradings of $\Hks$ and $H^{\frac{\infty}{2}+\gr}(A_0,M)$.\end{proposition}
\pr The proof of this proposition is parallel to that of Proposition \ref{hact}. In the proof we shall use the notion of the (not classical) derived functor of the functor of semiproduct. This derived functor is defined as follows.

Let $M\in \Armb$ be a right $A$--module.
First consider a functor $F_0:\oppAlmb \ra {\rm Vect}_\k$ defined on objects by
$$
F_0(M^\prime)=M\spr M^\prime,~~M^\prime \in \oppAlmb.
$$
This functor naturally extends to a functor $F:\KoppAlb \ra D({\rm Vect}_\k)$.

Fix an equivalence $\Phi:\DoppAlb \ra K({\mathcal S \mathcal J}(\oppA)_0)$ of the derived category $\DoppAlb$ with the homotopy category $K({\mathcal S \mathcal J}(\oppA)_0)$ of semijective complexes over $\oppAlmb$ provided by Theorem \ref{mainsinf}. From Corollary \ref{acycl} it follows that semijective complexes form a class of adapted objects for the functor $F$, and hence one can define the derived functor $DF:\DoppAlb \ra D({\rm Vect}_\k)$ of the functor $F$ as the composition  $DF=F\circ \Phi$. The functor $DF$ is exact and does not depend on the choice of the equivalence $\Phi$ (see Theorem 3.6 in \cite{V} for details).

Now observe that using Theorem \ref{threetor} and the standard semijective resolution $\sBaropp(\oppAo,N_0,\ke)$ of the $A_0$--module $\ke$ (see Proposition \ref{res2}) the semi--infinite cohomology space $H^{\frac{\infty}{2}+\gr}(A_0,M)$ may be calculated as follows:
$$
H^{\frac{\infty}{2}+\gr}(A_0,M)=H^\gr(M\spro \sBaropp(\oppAo,N_0,\ke)).
$$

From Proposition \ref{propspr} we also obtain that the complex $M\spro \sBaropp(\oppAo,N_0,\ke)$ for calculation of the semi--infinite cohomology space $H^{\frac{\infty}{2}+\gr}(A_0,M)$ may be represented as
$$
M\spr S_A\spro \sBaropp(\oppAo,N_0,\ke)=F(S_A\spro \sBaropp(\oppAo,N_0,\ke)).
$$
We shall prove that the complex of left $\oppA$--modules 
\begin{equation}
S_A\spro \sBaropp(\oppAo,N_0,\ke)
\end{equation}
 is semijective, and hence
\begin{equation}\label{beee}
H^{\frac{\infty}{2}+\gr}(A_0,M)=H^\gr(DF(S_A\spro \sBaropp(\oppAo,N_0,\ke))).
\end{equation}

Indeed, using the definition of the standard resolution $\sBaropp(\oppAo,N_0,\ke)$ and Lemma \ref{sbars} we have the following isomorphisms of complexes:
\begin{eqnarray*}
S_A\spro \sBaropp(\oppAo,N_0,\ke)=S_A\spro \sBar(A_0,N_0,S_{A_0}) \spro \ke = \\
\sBar(A_0,N_0,S_A)\spro \ke.
\end{eqnarray*}
Similarly to parts (a) and (b) of the proof of Proposition \ref{res2} (see Appendix 4) one can show that the complex $\sBar(A_0,N_0,S_A)\spro \ke$ is semijective. 
 
Now substituting the standard resolution $\sBaropp(\oppAo,N_0,\ke)$ of the one--dimen--sional $A_0$--module $\ke$ into the definition (\ref{sheckedef}) of Hecke algebras we obtain that
\begin{eqnarray*}
\Hks=\hDoppAl(S_A\spro \sBaropp(\oppAo,N_0,\ke)), \qquad \qquad  \\
\qquad \qquad \qquad  S_A\spro \sBaropp(\oppAo,N_0,\ke)).
\end{eqnarray*}

The algebra defined by the r.h.s. of the last equality naturally acts on the space defined by the r.h.s. of formula (\ref{beee}). Clearly, this action respects the gradings of $\Hks$ and $H^{\frac{\infty}{2}+\gr}(A_0,M)$. This completes the proof.

\qed


\subsection{W--algebras associated to complex semisimple Lie algebras as semi--infinite Hecke algebras}\label{Hklie}

\setcounter{equation}{0}
\setcounter{theorem}{0}

In this section we describe the W--algebras associated to complex semisimple Lie algebras as Hecke algebras. 
More precisely we shall identify the quantum BRST complex proposed in \cite{FF} for calculation of W--algebras with a standard complex for calculation of a semi--infinite Hecke algebra. 
This allows to obtain an invariant closed description of W--algebras without the bosonization technique used in \cite{FF}.

First we recall the definition of W--algebras associated to complex semisimple Lie algebras (see \cite{FF}, Sect. 4).
Let $\g$ be a complex semisimple Lie algebra,
${\ag}={\g}((z))\stackrel {\cdot}{+}{\mathbb C}K$ 
the non--twisted affine Lie algebra corresponding to $\g$. Recall that $\ag$ is the central extension of the loop algebra $\g((z))$ with the help of the standard two--cocycle $\omega_{st}$,
$$
\omega_{st}(x(z),y(z))={\rm Res}\langle x(z),y(z)\rangle\frac{dz}{z},
$$
where $\langle\cdot,\cdot\rangle$ is the Killing form of the Lie algebra $\g$.
 
Let $\n\subset \g$ be a maximal nilpotent subalgebra in $\g$ and $\an ={\n}((z))$ the loop algebra 
of the nilpotent Lie subalgebra ${\n}$. Note that $\an\subset {\ag}$ is a Lie subalgebra in $\ag$ because the standard cocycle $\omega_{st}$ vanishes when restricted to the subalgebra $\an ={\n}((z))\subset {\g}((z))$.  We denote by $\Uag$ and $\Uan$ the universal enveloping algebras of $\ag$ and $\an$, respectively. 

Let $\chi$ be the character of $\n$ which takes value $1$ on
all simple root vectors of $\n$ (see Example 1 in Section \ref{ex}). $\chi$ has a unique extension
to a character $\widehat \chi$ of  $\an ={\n}((z))$,
such that $\widehat \chi$ vanishes on the complement
$z^{-1} {\n}[[z^{-1}]] + z {\n}[[z]]$
of $\n$ in  ${\n}((z))$.
We denote by ${\mathbb C}_{\widehat \chi}$ the left one--dimensional $\Uan$--module that corresponds to $\widehat \chi$.

Let $\Uag_k$ be the quotient of the algebra $\Uag$ by the two--sided ideal generated by $K-k,~k\in {\mathbb C}$. Note that for any $k\in {\mathbb C}$ $\Uan$ is a subalgebra in $\Uag_k$ because the standard cocycle $\omega_{st}$ vanishes when restricted to the subalgebra $\an\subset {\g}((z))$.

Next observe that the algebras $\Uag_k$ and $\Uan$ inherit $\mathbb Z$--gradings from the natural $\mathbb Z$--gradings of $\ag$ and $\an$ by degrees of the parameter $z$, and satisfy conditions 
(i)--(viii) of Sections \ref{setup}, \ref{bimod}, with the natural triangular decompositions $\Uag_k=\Uagp_k\otimes \Uagm_k$ and $\Uan=\Uanp\otimes \Uanm$ provided by the decompositions $\ag= \ag_-+\ag_+$, $\an= \an_-+\an_+$, where $\ag_-=\g[[z^{-1}]]+{\mathbb C}K$, $\ag_+=z\g[[z]]$, $\an_{\pm}=\an\cap \ag_{\pm}$ (compare with Section \ref{lie}).  Hence one can define the algebras $\oppUag_k$, $\oppUan$ and  the semi--infinite Tor functors for $\Uag_k$ and $\Uan$. 

The algebra $\oppUag$ is explicitly described in the following proposition.

\begin{proposition}{\bf (\cite{Arkh1}, Proposition 4.6.7)}\label{oppUag}
The algebra $\oppUag_k$ is isomorphic to $\Uag_{-2h^\vee-k}$, where $h^\vee$ is the dual Coxeter number of $\g$.
\end{proposition}

Note also that from the explicit formula for the critical cocycle (see, for instance \cite{Arkh1}, Proposition 4.6.7) it follows that the critical cocycle of the algebra $\an$ is equal to zero. Therefore the trivial one--dimensional representation $\beta: \an \ra {\mathbb C}$ equips the algebra $\an$ with a semi--infinite structure (see Section \ref{lie}), and from Proposition \ref{oppug} we obtain that the algebra $\oppUan$ is isomorphic to $\Uan$.
We shall always identify the algebra $\oppUan$ with $\Uan$. 

Now consider the differential graded algebra
\begin{equation}\label{ff}
{\rm hom}_{\oppUag_k}^\gr((S_{\Uag_k}\underset{{\rm mod}-\an}{\otimes} {\mathbb C}_{\widehat \chi}^*)\otimes \sloo(\an),(S_{\Uag_k}\underset{{\rm mod}-\an}{\otimes} {\mathbb C}_{\widehat \chi}^*)\otimes \sloo(\an)),
\end{equation}
where $S_{\Uag_k}\underset{{\rm mod}-\an}{\otimes} {\mathbb C}_{\widehat \chi}^*$ is the tensor product in 
the category of right $\Uan$--modules of the right $\Uag_k$--module $S_{\Uag_k}$ and of the right 
one--dimensional $\Uan$--module ${\mathbb C}_{\widehat \chi}^*$ induced by the 
character $-\widehat \chi:\an \ra {\mathbb C}$; 
$(S_{\Uag_k}\underset{{\rm mod}-\an}{\otimes} {\mathbb C}_{\widehat \chi}^*)\otimes \sloo(\an)$ is the 
Feigin's standard complex for calculation of the semi--infinite cohomology of the right $\Uan$--module 
$S_{\Uag_k}\underset{{\rm mod}-\an}{\otimes} {\mathbb C}_{\widehat \chi}^*$ (see Corollary \ref{feigin}), 
and the $\oppUag_k$--module structure on this complex is induced by the left semiregular action of 
$\oppUag_k$ on $S_{\oppUag_k}$.

We shall show that the differential graded algebra (\ref{ff}) coincides with the quantum BRST 
complex, with the opposite multiplication, proposed in \cite{FF}, Sect. 4 for calculation of the 
W--algebra associated to the complex semisimple Lie algebra $\g$.
Indeed, consider the associative algebra
\begin{equation}\label{FFdga}
\Uag_k^{opp}\otimes C(\an +\an^*),
\end{equation}
where $C(\an +\an^*)$ is the Clifford algebra of the vector space $\an +\an^*$ equipped with 
the natural symmetric bilinear form (see Section \ref{lie}).

From the definitions of the semiregular bimodule $S_{\Uag_k}$ 
(see Section \ref{bimod}) 
and of the space $\sloo(\an)$ (see Section \ref{lie}) it follows that the 
algebra (\ref{ff}) coincides with the completion $\Uag_k^{opp}~\widehat \otimes ~C(\an +\an^*)$ 
of the algebra (\ref{FFdga})
by infinite series which are well defined as operators on the space 
$(S_{\Uag_k}\underset{{\rm mod}-\an}{\otimes} {\mathbb C}_{\widehat \chi}^*)\otimes \sloo(\an)$.
Here the action of the algebra $\Uag_k^{opp}$ on this space is induced by the 
right action of the algebra $\Uag_k$ on the space
$S_{\Uag_k}\underset{{\rm mod}-\an}{\otimes} {\mathbb C}_{\widehat \chi}^*$, and the Clifford algebra 
$C(\an +\an^*)$ naturally acts on the space $\sloo(\an)$ of semi--infinite exterior forms on $\g^*$.

Using this isomorphism we can equip the space $\Uag_k^{opp}~\widehat \otimes ~C(\an +\an^*)$ with the 
structure of a differential graded algebra. 
The differential in this algebra may be explicitly described as follows.

First the differential $d$ of the Feigin's standard complex 
$(S_{\Uag_k}\underset{{\rm mod}-\an}{\otimes} {\mathbb C}_{\widehat \chi}^*)\otimes \sloo(\an)$
(see formula (\ref{different}) in Section \ref{lie}) may be regarded as an element of
the graded algebra $\Uag_k^{opp}~\widehat \otimes ~C(\an +\an^*)$ of degree 1, the $\mathbb Z$--grading of the algebra $\Uag_k^{opp}~\widehat \otimes ~C(\an +\an^*)$ being induced by that of the differential graded algebra (\ref{ff}). Now from the definition
of the differential of the complex (\ref{ff}) (see formula (\ref{ddd})) it follows that the differential
of the graded algebra $\Uag_k^{opp}~\widehat \otimes ~C(\an +\an^*)$ induced by that of the 
differential graded algebra (\ref{ff}) is given by the supercommutator by element $d$.   
 
Using the last observation we conclude that the differential graded algebra 
$$
\Uag_k^{opp}~\widehat \otimes ~C(\an +\an^*)
$$ 
coincides with 
the quantum BRST complex defined in \cite{FF}, Sect. 4. 

By definition the W--algebra $W_k(\g)$ associated to the complex semisimple Lie algebra $\g$ is the opposite algebra of the
zeroth cohomology of the differential graded algebra (\ref{ff}),
\begin{eqnarray*}
W_k(\g)^{opp}= \qquad \qquad \qquad \qquad \qquad \qquad \qquad \qquad \qquad \qquad \qquad \qquad \qquad \qquad \qquad \\
\qquad H^0({\rm hom}_{\oppUag_k}^\gr((S_{\Uag_k}\underset{{\rm mod}-\an}{\otimes} {\mathbb C}_{\widehat \chi}^*)\otimes \sloo(\an),(S_{\Uag_k}\underset{{\rm mod}-\an}{\otimes} {\mathbb C}_{\widehat \chi}^*)\otimes \sloo(\an))).
\end{eqnarray*}

Now we realize the algebra $W_k(\g)$ as a semi--infinite Hecke algebra. First observe that the algebra $\Uag_k$ and the graded subalgebra $\Uan \subset \Uag_k$ satisfy the compatibility conditions of Section \ref{sheckedefr} under which the semi--infinite Hecke algebra of the triple $(\Uag_k,\Uan, {\mathbb C}_{\widehat \chi})$ may be defined.

\begin{proposition}\label{w}
The algebra $W_k(\g)^{opp}$ is isomorphic to the zeroth graded component of the semi--infinite Hecke algebra of the triple $(\Uag_k,\Uan, {\mathbb C}_{\widehat \chi})$,
$$
W_k(\g)^{opp}={\rm Hk}^{\frac{\infty}{2}+0}(\Uag_k,\Uan, {\mathbb C}_{\widehat \chi}).
$$
\end{proposition}
\pr
First we construct a standard complex for calculation of the Hecke algebra $\Hkso(\Uag_k,\Uan, {\mathbb C}_{\widehat \chi})$. 
Consider the standard semijective resolution (see Proposition \ref{stresU}) of the trivial one--dimensional $\Uan$--module,
$$
S_{\Uan}\otimes \sloo(\an).
$$
Then by Theorem 3.4 in \cite{V} the complex ${\mathbb C}_{\widehat \chi}\underset{\an-{\rm mod}}{\otimes} (S_{\Uan}\otimes \sloo(\an))$, where $\underset{\an-{\rm mod}}{\otimes}$ denotes the tensor product in the category of left $\Uan$--modules, is a semijective resolution of the  one--dimensional left $\Uan$--module ${\mathbb C}_{\widehat \chi}$.
Therefore by the definition of the semi--infinite Hecke algebra we have
\begin{eqnarray}\label{Hlie}
\Hkso(\Uag_k,\Uan, {\mathbb C}_{\widehat \chi})= \qquad \qquad \qquad \qquad \qquad \qquad \qquad \qquad \qquad  \\
\qquad \qquad \hDoppUagl(S_{\Uag_k}\spran({\mathbb C}_{\widehat \chi}\underset{\an-{\rm mod}}{\otimes} (S_{\Uan}\otimes \sloo(\an))), \nonumber \\
\qquad \qquad \qquad \qquad  \qquad \qquad \qquad S_{\Uag_k}\spran({\mathbb C}_{\widehat \chi}\underset{\an-{\rm mod}}{\otimes} (S_{\Uan}\otimes \sloo(\an)))). \nonumber
\end{eqnarray}

Similarly to Proposition \ref{stresU} one can show that the complex \\ $S_{\Uag_k}\spran({\mathbb C}_{\widehat \chi}\underset{\an-{\rm mod}}{\otimes}(S_{\Uan}\otimes \sloo(\an)))$ is semijective, and hence using Theorem \ref{mainsinf} we have
\begin{eqnarray*}
\Hkso(\Uag_k,\Uan, {\mathbb C}_{\widehat \chi})= \qquad \qquad \qquad \qquad \qquad \qquad \qquad \qquad \qquad \\
\qquad \qquad \hKoppUagl(S_{\Uag_k}\spran({\mathbb C}_{\widehat \chi}\underset{\an-{\rm mod}}{\otimes} (S_{\Uan}\otimes \sloo(\an))), \nonumber \\
\qquad \qquad \qquad \qquad  \qquad \qquad \qquad  S_{\Uag_k}\spran({\mathbb C}_{\widehat \chi}\underset{\an-{\rm mod}}{\otimes} (S_{\Uan}\otimes \sloo(\an)))). \nonumber
\end{eqnarray*}

Finally recalling formula (\ref{homK1}) 
the semi--infinite Hecke algebra of the triple $(\Uag_k,\Uan, {\mathbb C}_{\widehat \chi})$ may be calculated as the cohomology of the differential graded algebra
\begin{eqnarray}\label{ff1}
{\rm hom}_{\oppUag_k}^\gr(S_{\Uag_k}\spran({\mathbb C}_{\widehat \chi}\underset{\an-{\rm mod}}{\otimes} (S_{\Uan}\otimes \sloo(\an))), \\
\qquad \qquad \qquad S_{\Uag_k}\spran({\mathbb C}_{\widehat \chi}\underset{\an-{\rm mod}}{\otimes} (S_{\Uan}\otimes \sloo(\an)))). \nonumber
\end{eqnarray}
 
We shall show that there exists an isomorphism of complexes  of left $\oppUag_k$--modules,
\begin{eqnarray}\label{by}
S_{\Uag_k}\spran({\mathbb C}_{\widehat \chi}\underset{\an-{\rm mod}}{\otimes} (S_{\Uan}\otimes \sloo(\an)))=  \qquad \qquad \\
\qquad \qquad \qquad \qquad \qquad \qquad  (S_{\Uag_k}\underset{{\rm mod}-\an}{\otimes} {\mathbb C}_{\widehat \chi}^*)\otimes \sloo(\an), \nonumber
\end{eqnarray}
where $(S_{\Uag_k}\underset{{\rm mod}-\an}{\otimes} {\mathbb C}_{\widehat \chi}^*)\otimes \sloo(\an)$ is the Feigin's standard complex for calculation of semi--infinite cohomology of the right $\Uan$--module $S_{\Uag_k}\underset{{\rm mod}-\an}{\otimes} {\mathbb C}_{\widehat \chi}^*$.
The isomorphism (\ref{by}) provides an isomorphism of differential graded algebras (\ref{ff}) and (\ref{ff1}). As a consequence we obtain that the algebras $W_k(\g)^{opp}$ and ${\rm Hk}^{\frac{\infty}{2}+0}(\Uag_k,\Uan, {\mathbb C}_{\widehat \chi})$ are isomorphic.

In order to establish isomorphism (\ref{by}) it suffices to prove that
there exists an isomorphism of $\oppUag_k$--$\Uan$--bimodules,
\begin{equation}\label{biiso}
S_{\Uag_k}\spran({\mathbb C}_{\widehat \chi}\underset{\an-{\rm mod}}{\otimes} S_{\Uan})=S_{\Uag_k}\underset{{\rm mod}-\an}{\otimes} {\mathbb C}_{\widehat \chi}^*.
\end{equation}
This isomorphism is established similarly to the isomorphisms of Lemma \ref{propspr}.

First we recall two lemmas about modules over Lie algebras.
Let $M\in (\Uanp-{\rm mod})_0$ be a left $\Uanp$--module. Consider the tensor product of left $\Uanp$--modules $M$ and $\Uanp^*={\rm hom}_\k(\Uanp,\k)$  in the category of left $\Uanp$--modules, $M\underset{\an_+-{\rm mod}}{\otimes}\Uanp^*$. Here $\Uanp^*$ is equipped with the natural left action of $\Uanp$,
$$
u\cdot f(v)=f(vu),~~u,v\in \Uanp,~f\in \Uanp^*.
$$

Note that by Lemma \ref{tenshom} 
\begin{equation}\label{tensact}
M\underset{\an_+-{\rm mod}}{\otimes}\Uanp^*=M\otimes \Uanp^*= {\rm hom}_\k(\Uanp,M)
\end{equation}
as a vector space.
Equip the space ${\rm hom}_\k(\Uanp,M)$ with another left $\Uanp$--module structure as follows
\begin{equation}\label{oppact}
u\cdot f(v)=f(u^\top v),~~u,v\in \Uanp,~f\in {\rm hom}_\k(\Uanp,M).
\end{equation}
Here $\top$ stands for the canonical antiinvolution in $\Uanp$ defined on generators $x\in \an_+$ by $x^\top=-x$. We denote this left $\Uanp$--module by $M\otimes {\Uanp^{opp}}^*$.
 
\begin{lemma}{\bf (\cite{V}, Lemma 3.4)}\label{iso1}
Let $M\in (\Uanp-{\rm mod})_0$ be a left $\Uanp$--module. Then the map
$$
\begin{array}{l}
\phi: {\rm hom}_\k(\Uanp,M)\ra {\rm hom}_\k(\Uanp,M),\\
\\
\phi(f)(v)=\sum_i{v_1^i}^\top f({v_2^i}^\top),
\end{array}
$$
where $\sum_i{v_1^i}\otimes{v_2^i}=\Delta(v)$, and $\Delta$ is the comultiplication in $\Uanp$, $\Delta:\Uanp \ra \Uanp \otimes \Uanp$, provides both an isomorphism of left $\Uanp$--modules 
$$
M\underset{\an_+-{\rm mod}}{\otimes}\Uanp^* \ra M\otimes {\Uanp^{opp}}^*
$$
and the inverse isomorphism.
\end{lemma} 

Now let $L\in (\Uanm-{\rm mod})_0$ be a left $\Uanm$--module. Consider the tensor product of left $\Uanm$--modules $L$ and $\Uanm$  in the category of left $\an_-$--modules, $L\underset{\an_--{\rm mod}}{\otimes}\Uanm$. Here $\Uanm$ is equipped with the left regular action of $\Uanm$. 

Note that as a vector space  
\begin{equation}\label{tensact1}
L\underset{\an_--{\rm mod}}{\otimes}\Uanm=L\otimes \Uanm
\end{equation}
Equip the space $L\otimes \Uanm$ with another left module structure as follows
\begin{equation}\label{oppact1}
u\cdot (l\otimes v)=l\otimes vu^\top,~~u,v\in \Uanm,~l\in L.
\end{equation}
Here $\top$ stands for the canonical antiinvolution in $\Uanm$. We denote this left $\Uanm$--module by $L\otimes \Uanm^{opp}$
 
\begin{lemma}{\bf (\cite{V}, Lemma 3.2)}\label{iso2}
Let $L\in (\Uanm-{\rm mod})_0$ be a left $\Uanm$--module. Then the map
$$
\begin{array}{l}
\psi:L\otimes \Uanm \ra L\otimes \Uanm ,\\
\\
\psi(l\otimes v)=\sum_i{v_1^i}^\top l\otimes {v_2^i}^\top,
\end{array}
$$
where $\sum_i{v_1^i}\otimes{v_2^i}=\Delta(v)$, and $\Delta$ is the comultiplication in $\Uanm$, $\Delta:\Uanm \ra \Uanm \otimes \Uanm$, provides both an isomorphism of left $\Uanm$--modules 
$$
L\underset{\an_--{\rm mod}}{\otimes}\Uanm \ra L\otimes \Uanm^{opp}
$$
and the inverse isomorphism.
\end{lemma} 

Now we turn to the proof of isomorphism (\ref{biiso}).
First we calculate the space $S_{\Uag_k}\otimes^{\Uanp}({\mathbb C}_{\widehat \chi}\underset{\an-{\rm mod}}{\otimes} S_{\Uan})$. Using realization (\ref{oppSA2}) of the semiregular bimodule $S_{\Uan}$ and Lemma \ref{tenshom} we have the following isomorphisms of $\oppUag_k$--$\Uanm$--bimodules:
$$
\begin{array}{l}
S_{\Uag_k}\otimes^{\Uanp}({\mathbb C}_{\widehat \chi}\underset{\an-{\rm mod}}{\otimes} S_{\Uan})= \\
\\ 
S_{\Uag_k}\otimes^{\Uanp}({\mathbb C}_{\widehat \chi}\underset{\an-{\rm mod}}{\otimes}{\rm hom}_{\Uanm}(\Uan,\Uanm))= \\
\\
S_{\Uag_k}\otimes^{\Uanp}({\mathbb C}_{\widehat \chi}\underset{\an-{\rm mod}}{\otimes}{\rm hom}_{\k}(\Uanp,\Uanm))= \\
\\
S_{\Uag_k}\otimes^{\Uanp}({\mathbb C}_{\widehat \chi}\underset{\an-{\rm mod}}{\otimes}{\rm hom}_{\k}(\Uanp,\k))\otimes \Uanm.
\end{array}
$$

By Lemma \ref{iso1} the left $\Uanp$--module ${\mathbb C}_{\widehat \chi}\underset{\an-{\rm mod}}{\otimes}{\rm hom}_{\k}(\Uanp,\k)$ is isomorphic to the left $\Uanp$--module ${\mathbb C}_{\widehat \chi}\otimes {\Uanp^{opp}}^*$. 
Using this observation, the definition of the operation $\otimes^{\Uanp}$ and Lemma \ref{tenshom} we also have the following isomorphisms of left $\oppUag_k$--modules:
$$
\begin{array}{l}
S_{\Uag_k}\otimes^{\Uanp}({\mathbb C}_{\widehat \chi}\underset{\an-{\rm mod}}{\otimes}{\rm hom}_{\k}(\Uanp,\k))\otimes \Uanm= \\
\\
S_{\Uag_k}\otimes^{\Uanp}({\mathbb C}_{\widehat \chi}\otimes {\Uanp^{opp}}^*)\otimes \Uanm= \\
\\
{\rm hom}_{\Uanp}(\Uanp^{opp},S_{\Uag_k})\otimes {\mathbb C}_{\widehat \chi}\otimes \Uanm.
\end{array} 
$$

Finally the restriction isomorphism
$$
{\rm hom}_{\Uanp}(\Uanp^{opp},S_{\Uag_k})= S_{\Uag_k}
$$
yields an isomorphism of left $\oppUag_k$--modules, 
\begin{equation}\label{step11}
S_{\Uag_k}\otimes^{\Uanp}({\mathbb C}_{\widehat \chi}\underset{\an-{\rm mod}}{\otimes} S_{\Uan})= S_{\Uag_k}\otimes {\mathbb C}_{\widehat \chi}\otimes \Uanm
\end{equation}

Similarly, using realization (\ref{oppSA1}) of the semiregular bimodule $S_{\Uan}$ and Lemma \ref{iso2} we obtain an isomorphism of left $\oppUag_k$--modules,
\begin{equation}\label{step22}
S_{\Uag_k}\otimes_{\Uanm}({\mathbb C}_{\widehat \chi}\underset{\an-{\rm mod}}{\otimes} S_{\Uan})= S_{\Uag_k}\otimes {\mathbb C}_{\widehat \chi}\otimes \Uanp^*.
\end{equation}

Combining (\ref{step11}) and (\ref{step22}) and recalling the definition of the operation $\spran$ we conclude (see the proof of the same fact in Lemma \ref{propspr} for details) that
\begin{equation}\label{isom}
S_{\Uag_k}\spran{\mathbb C}_{\widehat \chi}\underset{\an-{\rm mod}}{\otimes} S_{\Uan}= S_{\Uag_k}\otimes {\mathbb C}_{\widehat \chi}.
\end{equation}
as a left $\oppUag_k$--module.

Using Lemmas \ref{iso1} and \ref{iso2} one checks directly that the right $\Uan$--action on the space $S_{\Uag_k}\spran({\mathbb C}_{\widehat \chi}\underset{\an-{\rm mod}}{\otimes} S_{\Uan})$ is transformed under isomorphism (\ref{isom}) into the action of $\Uan$ on the right $\Uan$--module
$$
S_{\Uag_k}\underset{{\rm mod}-\an}{\otimes} {\mathbb C}_{\widehat \chi}^*.
$$
This completes the proof of Proposition \ref{w}.

\qed

Using Proposition \ref{w} we shall explicitly calculate the algebra $W_k(\g)$.
Our main result in this section is
\begin{theorem}\label{mw}
The algebra $W_k(\g)^{opp}$ is canonically isomorphic to \\ ${\rm hom}_{\oppUag_k}(S_{\Uag_k}\spran {\mathbb C}_{\widehat \chi},S_{\Uag_k}\spran {\mathbb C}_{\widehat \chi})$,
\begin{equation}\label{wiso}
W_k(\g)^{opp}= {\rm hom}_{\oppUag_k}(S_{\Uag_k}\spran {\mathbb C}_{\widehat \chi},S_{\Uag_k}\spran {\mathbb C}_{\widehat \chi}).
\end{equation}
\end{theorem}

\begin{remark}
Recall that at the critical value of the parameter $k$, $k=-h^\vee$, the restricted completion $\cUag_{-h^\vee}$ of the algebra $\Uag_{-h^\vee}$ has a large center. This center is canonically isomorphic to the W-algebra $W_{-h^\vee}(\g)$ (see \cite{FF}, Proposition 6), 
$$
Z(\cUag_{-h^\vee})= W_{-h^\vee}(\g).
$$

From Theorem \ref{mw} we obtain a canonical algebraic isomorphism
\begin{equation}\label{zag}
Z(\cUag_{-h^\vee})= {\rm hom}_{\Uag_{-h^\vee}}(S_{\Uag_{-h^\vee}}\spran {\mathbb C}_{\widehat \chi},S_{\Uag_{-h^\vee}}\spran {\mathbb C}_{\widehat \chi})^{opp}. 
\end{equation}
Here using Proposition \ref{oppUag} we replaced the algebra $\oppUag_{-h^\vee}$ with ${\Uag_{-h^\vee}}$ (We note that at the critical level of the parameter $k$, $k=-h^\vee$ the algebra $\Uag_{-h^\vee}$ is ``self--dual'' in the sense that the algebra $\oppUag_{-h^\vee}$ is isomorphic to ${\Uag_{-h^\vee}}$).

The  description (\ref{zag}) of the center $Z(\cUag_{-h^\vee})$
is similar to the realization of the center $Z(\Ug)$ of the algebra $\Ug$ obtained by Kostant in \cite{K} (see Example 1 in Section \ref{ex}).
\end{remark}

\pr
First by Proposition \ref{w} the algebra $W_k(\g)^{opp}$ is isomorphic to the zeroth graded component ${\rm Hk}^{\frac{\infty}{2}+0}(\Uag_k,\Uan, {\mathbb C}_{\widehat \chi})$ of the semi--infinite Hecke algebra of the triple $(\Uag_k,\Uan, {\mathbb C}_{\widehat \chi})$.
We shall apply vanishing theorem \ref{vv} to calculate the algebra $\Hkso(\Uag_k,\Uan, {\mathbb C}_{\widehat \chi})$,
\begin{eqnarray*}
\Hkso(\Uag_k,\Uan, {\mathbb C}_{\widehat \chi})= \qquad \qquad \qquad \qquad \qquad \qquad \qquad  \\
\qquad \qquad \hDoppUagl(S_{\Uag_k}\spran S^\gr({\mathbb C}_{\widehat \chi}),  \\
\qquad \qquad \qquad \qquad \qquad  \qquad \qquad \qquad S_{\Uag_k}\spran S^\gr({\mathbb C}_{\widehat \chi})), 
\end{eqnarray*}
where $S^\gr({\mathbb C}_{\widehat \chi})$ is a semijective resolution of the left $\Uan$--module ${\mathbb C}_{\widehat \chi}$.

More precisely, we shall show that the right $\Uan$--module $S_{\Uag_k}$ is $\Uanp$--injective and $\Uanm$--projective. Then by Proposition \ref{vv} 
$$
H^\gr(S_{\Uag_k}\spran S^\gr({\mathbb C}_{\widehat \chi}))= S_{\Uag_k}\spran {\mathbb C}_{\widehat \chi},
$$
and hence Proposition \ref{svanish} yields an algebraic isomorphism,
\begin{eqnarray*}
\Hkso(\Uag_k,\Uan, {\mathbb C}_{\widehat \chi})= \qquad \qquad \qquad \qquad \qquad \qquad \qquad  \\
\qquad \qquad \hDoppUagl(S_{\Uag_k}\spran {\mathbb C}_{\widehat \chi},S_{\Uag_k}\spran {\mathbb C}_{\widehat \chi}). 
\end{eqnarray*}
In particular, this establishes isomorphism \ref{wiso}.

Now we prove that the right $\Uan$--module $S_{\Uag_k}$ is $\Uanp$--injective and $\Uanm$--projective.
Using realizations (\ref{SA1}) and (\ref{SA2}) of the right $\Uag_k$--module $S_{\Uag_k}$ we obtain that
\begin{equation}\label{r1}
S_{\Uag_k}=\Uagp_k^*\otimes \Uagm_k
\end{equation}
as a right $\Uagm_k$--module, and
\begin{equation}\label{r2}
S_{\Uag_k}={\rm hom}_\k(\Uagp_k,\Uagm_k)
\end{equation}
as a right $\Uagp_k$--module.

Now observe that the right $\Uagm_k$--module $\Uagm_k$ is $\Uanm$--projective because $\an_-$ is a Lie subalgebra in $\ag_-$ (see \cite{carteil}, Ch. XIII, Proposition 4.1). Therefore $\Uagp_k^*\otimes \Uagm_k$ is also projective as a right $\Uanm$--module.

 Similarly the right $\Uagp_k$--module $\Uagp_k$ is $\Uanp$--projective because $\an_+$ is a Lie subalgebra in $\ag_+$, and hence the the right $\Uagp_k$--module ${\rm hom}_\k(\Uagp_k,\Uagm_k)$ is $\Uanp$--injective (see \cite{carteil}, Ch. VI, Proposition 1.4).
This completes the proof of the theorem.

\qed


\section*{Open problems and concluding remarks}

In conclusion we discuss possible applications of the Hecke algebras defined in this paper and formulate some interesting open problems.

In this paper we only studied general properties of Hecke algebras. However in particular situations these properties may have nontrivial important consequences. For instance, it would be interesting to obtain a proof of the duality theorem for W--algebras (see \cite{FF}, Proposition 5) and the isomorphism of the algebra $W_{-h^\vee}(\g)$ with the corresponding classical Poisson  W--algebra (see \cite{FF}, Sect. 5) using the description of the W--algebras obtained in Theorem \ref{mw} of this paper. Note that the distinguished role of the critical level of the parameter $k$, $k=-h^\vee$ may be easily observed in this description (see Remark after Theorem \ref{mw}).

It also would be interesting to apply the construction of semi--infinite Hecke algebras to define opers which are further generalizations of W--algebras playing an important role in the Drinfeld's approach to the geometric Langlands correspondence.

Another possible important application of semi--infinite Hecke algebras is q--deformation of W--algebras. These results will be presented in a subsequent paper.

In this paper we do not discuss applications of Hecke algebras to quantum mechanics and quantum field theory. These applications are connected with quantum reduction procedure (see, for instance, \cite{D} or \cite{S}, Sect. 6). 

A series of open problems is concerned with semi--infinite cohomology of associative algebras. The first question is related to the technical conditions imposed on algebra $A$ to define the semi--infinite Tor functor (see conditions (i)--(viii) of Sections \ref{setup}, \ref{bimod}). The algebras satisfying these conditions play an important role not only in the semi--infinite cohomology theory (see, for instance, \cite{man}). But it seems that in the definition of the semi--infinite Tor functor these conditions may be made weaker.

Indeed, observe that the algebra $B$ entering the triangular decompositions for the algebra $A$ (see condition (vi) of Section \ref{setup}) does not explicitly appear in the definition of semijective complexes. Moreover, in the proof of the main theorem of semi--infinite homological algebra (Theorem \ref{mainsinf} in this paper) the triangular decompositions are only used to show that every module from the category $\Almb$ is a submodule of an $N$--injective module and a strong quotient of a relatively to $N$ projective module (see \cite{V}, proof of Lemma 3.6, step 1 and proof of Lemma 3.7, step 1). Therefore Theorem \ref{mainsinf} holds, in fact, under more weak conditions.

Moreover, in the definition of the right semiregular representation $S_A=N^*\otimes_NA$ (see Section \ref{bimod}) we only used the subalgebra $N$. Therefore the semiregular bimodule may be defined under weaker conditions as well.  However in this case properties of this bimodule are not clear.

In \cite{Arkh3} S. Arkhipov also tries to introduce a functor that corresponds to the semiproduct functor in our terminology (see Corollary 4.7.2 in \cite{Arkh3}). Remarkably, his definition only uses the subalgebra $N\subset A$. Unfortunately this definition does not exactly coincide with the definition of the semiproduct functor, but probably may be corrected.

From these observations we conclude that there is a possibility to define the semi--infinite Tor functor under conditions weaker than conditions (i)--(viii) of Sections \ref{setup}, \ref{bimod}, e.g. using subalgebra $N\subset A$ only.

Another series of questions is connected with properties of the semi--infinite cohomology. In this paper we only discussed the properties of the semi--infinite Tor functor which follow from the definition of this functor and from the basic properties of the semiregular bimodule. Essentially these properties are of homotopic nature.

But the vanishing theorem proved by Voronov in \cite{V} (see also Proposition \ref{vv} in this paper) shows that the semi--infinite Tor functor has some unexpected properties. The Voronov's vanishing theorem seems to be connected with properties of the semiproduct functor with respect to quasi--isomorphisms, i.e. this vanishing property is not of homotopic nature.

A related question is connected with the partial and the full (not classical) derived functors of the functor of semiproduct. As we observed in the remark after Theorem \ref{threetor} the corresponding classical derived functors are not isomorphic. But there is a possibility that the nonclassical partial and full derived functors of the functor of semiproduct are isomorphic. This property is not of homotopic nature even in case of the usual tensor product functor (see \cite{GM}, Ch. III, \S 7, Ex. 6).


\section*{Appendix}


\subsection*{A1. Two lemmas about ``almost'' double complexes}

\renewcommand{\thetheorem}{A1.\arabic{theorem}}
\renewcommand{\theequation}{A1.\arabic{equation}}
\renewcommand{\thelemma}{A1.\arabic{lemma}}

\setcounter{equation}{0}
\setcounter{theorem}{0}

In this section we prove two technical lemmas about ``almost'' double complexes which are ``partially'' bounded.

\begin{lemma}\footnote{A modification of this lemma was incorrectly formulated and proved in \cite{Arkh1}, Lemma 3.2.2.}\label{abovebound}
Let $X^{\gr}=\bigoplus_{n\in {\mathbb Z}}X^n,~~X^n=\prod_{p\in {\mathbb Z}}X^{-p,p+n}$ be a complex over an abelian category $\mathcal A$ with differential $d=d_1+d_2;~~d_1:X^{p,q}\ra X^{p+1,q},~d_2:X^{p,q}\ra X^{p,q+1}$. Suppose that for every $q\in {\mathbb Z}$ the complex $(X^{\gr , q},d_1)$ is homotopic to zero, and $X^{\gr ,q}=0$ for $q>N,~N\in {\mathbb Z}$. Then the complex $X^{\gr}$ is acyclic.
\end{lemma}
\pr 
Choose homotopy maps for the  complexes $(X^{\gr , q},d_1)$, 
$$
h_q:X^{p,q}\ra X^{p+1,q},~~h_qd_1+d_1h_q={\rm Id}_{X^{\gr , q}}.
$$  

Clearly, the sum of these maps $h=\sum_{q\in {\mathbb Z}}h_q$ is a well defined map $h:X^n\ra X^{n-1}$ such that
\begin{equation}\label{homzero}
dh+hd=Id_{X^{\gr}}+d_2h+hd_2.
\end{equation}

Since $X^{\gr ,q}=0$ for $q>N$ the map $d_2h+hd_2:X^\gr\ra X^\gr$ is nilpotent, and hence the map $Id_{X^{\gr}}+d_2h+hd_2$ is invertible being the sum of the identity map and of a nilpotent one. Therefore this map induces an invertible map of the cohomology space of the complex $X^\gr$.

On the other hand by (\ref{homzero}) the map $Id_{X^{\gr}}+d_2h+hd_2$ is homotopic to zero, and hence it induces the zero map of the cohomology space $H^\gr(X^\gr)$. 

Finally we obtain that the zero map of space $H^\gr(X^\gr)$ onto itself is invertible. We conclude that $H^\gr(X^\gr)=0$. This completes the proof.

\qed

\begin{lemma}\label{belowbound}
Let $X^{\gr}=\bigoplus_{n\in {\mathbb Z}}X^n,~~X^n=\prod_{p\in {\mathbb Z}}X^{-p,p+n}$ be a complex over an abelian category $\mathcal A$ with differential $d=d_1+d_2;~~d_1:X^{p,q}\ra X^{p+1,q},~d_2:X^{p,q}\ra X^{p,q+1}$. Suppose that for every $q\in {\mathbb Z}$ the complex $(X^{\gr , q},d_1)$ is acyclic, and the only non--vanishing components $X^{p,q}$ are situated in the upper--half of the p--q plane, i.e. $X^{\gr ,q}=0$ for $q<0$. Then the complex $X^{\gr}$ is acyclic.
\end{lemma}

\begin{remark}
If $X^{\gr, \gr}$ is a double complex then a similar result for this complex holds if it is situated in the lower--half of the p--q plane (compare with \cite{MacLane}, Ch. XI, \S 6, Theorem 6.1), i.e. replacing direct products by direct sums in the definition of graded components $X^n$ converts the grading condition of Lemma \ref{belowbound} to the opposite one.
\end{remark}
 
\pr
Let $u=\sum_{p=-\infty}^{p_0}u^{p,n-p},~u^{p,n-p}\in X^{p,n-p},~p=-\infty,\ldots,p_0$ be a cocycle in $X^\gr$, i.e. $du=0$. We shall construct an element $v=\sum_{p=-\infty}^{p_0-1}v^{p,n-1-p},~v^{p,n-1-p}\in X^{p,n-1-p},~p=-\infty,\ldots,p_0-1$ such that $u=dv$. In order to do that we shall use induction procedure based on diagram chasing in the following commutative diagram (Fig.1): 

$$
\xymatrix@!
{
*+[o][F-]{v}\ar@{--}[1,1]   & *+[o][F-]{u} \ar@{--}[1,1]  & {~~~~}   & {~~~~}  & {~~~~}  & {~~}  & {~~~~}  & {~~~~}  \\
 {~~~~}  & {\gr}\ar[r] \ar@{--}[u] \ar@{-} '[1,1] '[2,2] [4,4]^(0.5){n-1}  &  {\gr}\ar@{-->}[r]\ar@{--}[u] \ar@{-} '[1,1] '[2,2] [4,4]^(0.5){n}  & {0}  & {~~~~}  & {~~~~}  & {~~~~}  & {~~~~}   \\
 {~~~~}  & {~~~~}  &  {\gr}\ar[r] \ar[u]  &  {\gr}\ar@{-->}[r] \ar@{-->}[u]  & {0}  & {~~~~}  & {~~~~}  & {~~~~} \\
 {~~~~}  & {~~~~}  & {~~~~}  & {\gr} \ar[r] \ar[u]   & {\gr}\ar@{-->}[r] \ar@{-->}[u]  & {0}  & {~~~~}   & {~~~~}   \\
 {~~~~}\ar[0,7]_>*+{p,d_1}  & {~~~~}  & {~~~~}  & {~~~~}  & {~~~~}  & {~~~~}  & {~~~~}  & {~~~~}  & \\
 {~~~~}  & {~~~~}  & {~~~~}  & {~~~~}  & {~~~~}  & {~~~~}  & {~~~~}\ar[-5,0]_>*+{q,d_2} & {~~~~} 
  }
$$
\begin{center}
 Fig.1
\end{center}
Here the components of $u$ and $v$ are situated at the black dots on the lines marked by framed $u$ and $v$, respectively. The condition $du=0$ is equivalent to the fact that the sum of the images of every two dashed arrows having the same target in this diagram is equal to zero. 

Formally the induction procedure looks like as follows. 
First the equality $du=0$ being rewritten in components takes the form:
\begin{equation}\label{du0}
d_1u^{p_0,n-p_0}=0, 
\end{equation}
\begin{equation}\label{du}
d_1u^{p_0-k-1,n-p_0+k+1}+d_2u^{p_0-k,n-p_0+k}=0,~~k=0,\ldots,\infty .
\end{equation}

Since by (\ref{du0}) $d_1u^{p_0,n-p_0}=0$ and 
for every $q\in {\mathbb Z}$ the complex $(X^{\gr , q},d_1)$ is acyclic one can find $v^{p_0-1,n-p_0}\in X^{p_0-1,n-p_0}$ such that 
\begin{equation}\label{nol}
u^{p_0,n-p_0}=d_1v^{p_0-1,n-p_0}.
\end{equation}

Substituting (\ref{nol}) into (\ref{du}) for $k=0$ we obtain that
\begin{equation}\label{pass}
\begin{array}{l}
d_1u^{p_0-1,n-p_0+1}+d_2u^{p_0,n-p_0}= \\
\\
d_1u^{p_0-1,n-p_0+1}+d_2d_1v^{p_0-1,n-p_0}= \\
\\
d_1u^{p_0-1,n-p_0+1}-d_1d_2v^{p_0-1,n-p_0}= \\
\\
d_1(u^{p_0-1,n-p_0+1}-d_2v^{p_0-1,n-p_0})=0.
\end{array}
\end{equation}
Hence one can find an element $v^{p_0-2,n-p_0+1}\in X^{p_0-2,n-p_0+1}$ such that 
\begin{equation}\label{odin}
u^{p_0-1,n-p_0+1}-d_2v^{p_0-1,n-p_0}=d_1v^{p_0-2,n-p_0+1}.
\end{equation}

Summarizing (\ref{nol}) and (\ref{odin}) gives:
$$
u^{p_0,n-p_0}+u^{p_0-1,n-p_0+1}=dv^{p_0-1,n-p_0}+d_1v^{p_0-2,n-p_0+1}.
$$
This completes the base of the induction.

Now suppose that 
\begin{equation}\label{indass1}
\sum_{k=0}^{N}u^{p_0-k,n-p_0+k}=d(\sum_{k=0}^{N-1}v^{p_0-k-1,n-p_0+k})+d_1v^{p_0-N-1,n-p_0+N},
\end{equation}
where elements $v^{p_0-k-1,n-p_0+k}\in X^{p_0-k-1,n-p_0+k},~k=0,\ldots,N$ satisfying the following system of equations
\begin{equation}\label{indass}
u^{p_0-k,n-p_0+k}=d_2v^{p_0-k,n-p_0+k-1}+d_1v^{p_0-k-1,n-p_0+k},~~
k=1,\ldots,N 
\end{equation}
have already been constructed.

We shall find $v^{p_0-N-2,n-p_0+N+1}\in X^{p_0-N-2,n-p_0+N+1}$ such that
\begin{equation}\label{toprove}
\sum_{k=0}^{N+1}u^{p_0-k,n-p_0+k}=d(\sum_{k=0}^{N}v^{p_0-k-1,n-p_0+k})+d_1v^{p_0-N-2,n-p_0+N+1}. 
\end{equation} 

Using (\ref{du}) for $k=N-1$ and equality (\ref{indass}) for $k=N$ we obtain, similarly to the first step of the induction (see (\ref{pass})), that
$$
d_1(u^{p_0-N-1,n-p_0+N+1}-d_2v^{p_0-N-1,n-p_0+N})=0.
$$
Therefore one can find an element $v^{p_0-N-2,n-p_0+N+1}\in X^{p_0-N-2,n-p_0+N+1}$ such that 
\begin{equation}\label{dva}
u^{p_0-N-1,n-p_0+N+1}=d_2v^{p_0-N-1,n-p_0+N}+d_1v^{p_0-N-2,n-p_0+N+1}.
\end{equation}

Now adding (\ref{indass}) and (\ref{dva}) yields representation (\ref{toprove}) for the sum \\ $\sum_{k=0}^{N+1}u^{p_0-k,n-p_0+k}$.

Iterating this procedure we finally obtain that
$$
u=d(\sum_{p=-\infty}^{p_0-1}v^{p,n-1-p})
$$
for some $v^{p,n-1-p}\in X^{p,n-1-p}$. This concludes the proof.

\qed


\subsection*{A2. Special inverse and direct limits}

\renewcommand{\thetheorem}{A2.\arabic{theorem}}
\renewcommand{\theequation}{A2.\arabic{equation}}
\renewcommand{\thelemma}{A2.\arabic{lemma}}
\renewcommand{\thedefinition}{A2.\arabic{definition}}

\setcounter{equation}{0}
\setcounter{theorem}{0}
\setcounter{definition}{0}

In this section we recall, following \cite{Sp}, the notions of special inverse and direct limits.

Let $\mathcal A$ be an abelian category, $\mathcal J \subset {\rm Kom}({\mathcal A})$ a class of complexes.

\begin{definition}
An inverse system $(I^\gr_n)_{n\in {\mathbb Z}}$ in $\mathcal J$ is called a $\mathcal J$--special inverse system if for every $n\in {\mathbb Z}$ the natural map 
$$
I^\gr_n\ra I^\gr_{n-1}
$$
is surjective, its kernel $C^\gr_n$ belongs to $\mathcal J$ and the short exact sequence 
$$
0\ra C^\gr_n\ra I^\gr_n\ra I^\gr_{n-1}\ra 0
$$
is split in each degree.
\end{definition}

\begin{definition}
The class $\mathcal J$ is closed under special inverse limits if every $\mathcal J$--special inverse system in $\mathcal J$ has a limit which is contained in $\mathcal J$, and every complex isomorphic in ${\rm Kom}({\mathcal A})$ to a complex in $\mathcal J$ is contained in $\mathcal J$.
\end{definition}

The following proposition provides important examples of classes of complexes closed under special inverse limits
\begin{proposition}{\bf (\cite{Sp}, Corollary 2.5)}
Let $\mathcal I\subset \mathcal A$ be a class of complexes. Then the class of all complexes $S^\gr \in {\rm Kom}(\mathcal A)$ such that ${\rm Hom}_{K(\mathcal A)}^\gr(A^\gr, S^\gr)=0$ for every $A^\gr \in \mathcal I$ is closed under special inverse limits.
\end{proposition} 

Applying this proposition to the category $\mathcal A =(N-{\rm mod})_0$ and the class $\mathcal I$ of all acyclic complexes in ${\rm Kom}(N)_0$ we obtain
\begin{corollary}{\bf (\cite{V}, Lemma 3.6, step 4)}\label{spinv}
Let $A$ be a $\mathbb Z$-graded associative algebra over a field $\k$ containing a graded subalgebra $N$. Then the class of $N$--K--injective complexes in $\CAlb$ is closed under special inverse limits.
\end{corollary}

Dually one can define the notion of special direct limits.
Let $\mathcal P \subset {\rm Kom}({\mathcal A})$ be a class of complexes.

\begin{definition}
A direct system $(P^\gr_n)_{n\in {\mathbb Z}}$ in $\mathcal P$ is called a $\mathcal P$--special direct system if for every $n\in {\mathbb Z}$ the natural map 
$$
P^\gr_{n-1}\ra P^\gr_{n}
$$
is injective, its cokernel $C^\gr_n$ belongs to $\mathcal P$ and the short exact sequence 
$$
0\ra P^\gr_{n-1}\ra P^\gr_{n}\ra C^\gr_n\ra 0
$$
is split in each degree.
\end{definition}

\begin{definition}
The class $\mathcal P$ is closed under special direct limits if every $\mathcal P$--special direct system in $\mathcal P$ has a limit which is contained in $\mathcal P$, and every complex isomorphic in ${\rm Kom}({\mathcal A})$ to a complex in $\mathcal P$ is contained in $\mathcal P$.
\end{definition}

The following proposition provides important examples of classes of complexes closed under special direct limits
\begin{proposition}{\bf (\cite{Sp}, Corollary 2.8)}
Let $\mathcal I\subset \mathcal A$ be a class of complexes. Then the class of all complexes $S^\gr \in {\rm Kom}(\mathcal A)$ such that ${\rm Hom}_{K(\mathcal A)}^\gr(S^\gr, A^\gr)=0$ for every $A^\gr \in \mathcal I$ is closed under special direct limits.
\end{proposition} 

Applying this proposition to the category $\mathcal A=\Almb$ and the class $\mathcal I$ of all acyclic complexes in $\CAlb$ which are isomorphic to zero in the category $K(N)_0$ we obtain
\begin{corollary}{\bf (\cite{V}, Lemma 3.7, step 4)}\label{spdir}
Let $A$ be a ${\mathbb Z}$--graded associative algebra over a field $\k$ containing a graded subalgebra $N$. Then the class of relative to $N$ K--projective complexes in $\CAlb$ is closed under special direct limits.
\end{corollary}


\subsection*{A3. Proof of Proposition \ref{res1}}

\renewcommand{\theequation}{A3.\arabic{equation}}

\setcounter{equation}{0}

In the proof of Proposition \ref{res1} we shall use Lemmas \ref{abovebound} and \ref{belowbound} proved in the Appendix.

First we verify that the complex $\sBar(A,N,M)$ is K-injective as a complex of $N$--modules. By part (iii) of Proposition \ref{barprop} the complex $\Bar(A,N,A)$ is homotopically equivalent to $A$ as a complex of $A$--$N$--bimodules. It follows that the complex 
$$\sBar(A,N,M)=\hAd(\Bar(A,B,A),M)\otimes_A\Bar(A,N,A)
$$
is homotopically equivalent to $\hAd(\Bar(A,B,A),M)$ as a complex of $N$--modules. But from the definition of the standard normalized bar resolution (see Section \ref{stres}) and the triangular decompositions for the algebra $A$ (see condition (vi) of Section \ref{setup}) we also have a natural isomorphism of complexes of $N$--modules,
$$
\hAd(\Bar(A,B,A),M)={\rm hom}_N^{\gr}(\Bar(N,{\k},N),M).
$$

The complex ${\rm hom}_N^{\gr}(\Bar(N,{\k},N),M)$  is obviously a bounded from below complex of $N$--injective modules. By part 1 of Proposition \ref{sinjprop} this complex is also K-injective. This proves that the complex $\sBar(A,N,M)$ is K-injective as a complex of $N$--modules.

Next we show that $\sBar(A,N,M)$ is K--projective relative to $N$. By definition we have to prove that for every complex of right $A$--modules $V^\gr\in {\rm Kom}\Armb$, that is homotopically equivalent to zero as a complex of $N$--modules, the complex $\hAd(\sBar(A,N,M),V^\gr)$ is acyclic.
Denote this complex by $X^\gr$, 
$$
X^\gr=\HAd(\sBar(A,N,M),V^\gr).
$$

To apply Lemmas \ref{abovebound} and \ref{belowbound} to the complex $X^\gr$ we have to explicitly rewrite this complex in components. Using the definition of the complex $\sBar(A,N,M)$ we have:
\begin{eqnarray}\label{be}
X^\gr=\bigoplus_{q\in {\mathbb Z}}\prod_{l+p=q}\prod_{n+k=p}X^{l,k,n},
\end{eqnarray} 
where 
$$
X^{l,k,n}=\hA(\hA(\Baro^k(A,B,A),M)\otimes_A\Barn(A,N,A),V^l).
$$

Note that $X^\gr$ is an ``almost'' three--complex with components $X^{l,k,n}$. 
We denote by $d_1,~d_2$ and $d_3$ the differentials in $X^\gr$ induced by the differentials of the complexes $V^\gr,~\Bar(A,B,A)$ and $\Bar(A,N,A)$, respectively. The total differential $d$  of $X^\gr$ is the sum $d=d_1+d_2+d_3$.

First using Lemma \ref{abovebound} we calculate the cohomology of the complex $X^\gr$ with respect to the differential $d'=d_1+d_2$. In order to do that we note that by the definition of the normalized bar resolution $\Bar(A,N,A)$ (see formula (\ref{bar})) the complex (\ref{be}), with the differential $d_3$ forgotten, may be rewritten as 
\begin{eqnarray*}\label{star}
X^\gr=\bigoplus_{q\in {\mathbb Z}}\prod_{l+p=q}\prod_{n+k=p}{\rm hom}_N(\hA(\Baro^k(A,B,A),M)\otimes_A \qquad \\
\qquad \qquad T^{n+1}_0(A,N),V^l), \nonumber
\end{eqnarray*}
where $T^{n+1}_0(A,N)$ is the quotient of the tensor product 
$$
T^{n+1}(A,N)=\underbrace{A\otimes_N\ldots \otimes_NA}_{n+1 ~\mbox{\tiny  times}}
$$ 
by the subspace $\overline T^{n+1}(A,N)$,
$$
\overline T^{n+1}(A,N)=\{ a_0\otimes\ldots\otimes a_n\in T^{n+1}(A,N) |~\exists s\in\{1,\ldots,n\}: a_s\in N\}.
$$

Now recall that the complex $V^\gr$ is homotopically equivalent to zero as a complex of $N$--modules. Observe also that the functor ${\rm hom}_N$ preserves this property. Therefore for every $k,n\in {\mathbb Z}$ the complex $X^{\gr,k,n}$ is homotopic to zero. Note also that $X^{\gr,k,n}=0$ for $k>0$, and hence from Lemma \ref{abovebound}, applied to the complex $X^\gr$ equipped with the differential $d'=d_1+d_2$, it follows that $H^{\gr}(X^\gr,d')=0$.

Next, we apply Lemma \ref{belowbound} to the complex $(X^{\gr},d)$, the differential $d$ being split as follows $d=d'+d_3$. It suffices to remark that $X^{l,k,n}=0$ for $n<0$ and $H^{\gr}(X^\gr,d')=0$ as we proved above. Therefore by  Lemma \ref{belowbound} $H^\gr(X^\gr,d)=0$. This proves that the complex $\sBar(A,N,M)$ is K--projective relative to $N$. But we have also proved that this complex is K--injective as a complex of $N$--modules. We conclude that $\sBar(A,N,M)$ is a semijective complex. 

Now we prove that the semijective complex $\sBar(A,N,M)$ is a semijective resolution of $M$. By proposition \ref{resmod} it suffices to verify that $H^\gr(\sBar(A,N,M))=M$. 
Indeed, the complex $\hAd(\Bar(A,B,A),M)\otimes_A\Bar(A,N,A)$ is quasi--isomor--phic to $\hAd(\Bar(A,B,A),M)$ (see \cite{GM}, Lemma III.7.12), and the complex \\ $\hAd(\Bar(A,B,A),M)$ is quasi--isomorphic to $M$ by definition. The last two facts imply, in particular, that $H^\gr(\sBar(A,N,M))=M$. This completes the proof.

\qed


\subsection*{A4. Proof of Proposition \ref{res2}}

\renewcommand{\thetheorem}{A4.\arabic{theorem}}
\renewcommand{\theequation}{A4.\arabic{equation}}
\renewcommand{\thelemma}{A4.\arabic{lemma}}
\renewcommand{\thedefinition}{A4.\arabic{definition}}

\setcounter{equation}{0}
\setcounter{theorem}{0}
\setcounter{definition}{0}

The proof will be divided into three steps:

In part (a) we show that $\sBaropp (\oppA ,N,M^\prime)$ is K-injective as a complex of $N$--modules.  In part (b) we prove that this complex is also relatively to $N$ K--projective. In part (c) we verify that $H^\gr(\sBaropp (\oppA ,N,M^\prime))=M^\prime$. By Proposition \ref{resmod} properties (a), (b) and (c) imply that $\sBaropp (\oppA ,N,M^\prime)$ is a semijective resolution of $M^\prime$.

(a) First we prove that the complex 
$$
\sBaropp(\oppA,N,M^\prime)=\sBar(A,N,S_A)\spr M^\prime
$$ 
is K-injective as a complex of $N$--modules. We construct a special inverse system of $N$--K--injective complexes converging to $\sBaropp(\oppA,N,M^\prime)$ (see \cite{Sp} or Appendix 2 in this paper for the definition of special inverse systems). By  Corollary \ref{spinv} this implies that the complex $\sBaropp(\oppA,N,M^\prime)$ is  $N$--K--injective itself.

The required special inverse system is defined as follows. Observe that \\ $\sBaropp(\oppA,N,M^\prime)$ is the total complex of the bicomplex
$$
K^{p,q}=\hA(\Baro^{-p}(A,B,A),S_A)\otimes_A\Baro^q(A,N,A)\spr M^\prime,~p\geq 0,~q\leq 0.
$$
We denote by $d_1$ and $d_2$ the differentials in $K^{p,q}$ induced by the differentials of $\Bar(A,B,A)$ and $\Bar(A,N,A)$, respectively.  

Consider a set $I^\gr_{n},~n\in {\mathbb Z}$ of complexes in ${\rm Kom}(\oppA)_0$, where $I^\gr_{n}$ is the total complex of the bicomplex $\tau^2_{\geq -n}K^{\gr,\gr}$, and $\tau^2_{\geq -n}K^{\gr,\gr}$ denotes the truncation from below of the bicomplex  $K^{p,q}$ with respect to the second grading,
$$
\tau^2_{\geq -n}K^{\gr,\gr}=0\ra {\rm Coker}~d_2^{-n-1}\ra K^{\gr,-n+1}\ra \ldots .
$$
Clearly, $I^\gr_{n},~n\in{\mathbb Z}$ is an inverse system of complexes converging to $\sBaropp(\oppA,N,M^\prime)$. We also note that the natural maps $I^\gr_n\ra I^\gr_{n-1}$ are surjective, their kernels are given by:
$$
C_n^\gr=0\ra {\rm Coker}~d_2^{-n-1}\ra {\rm Im}~d_2^{-n}\ra 0.
$$ 
Therefore in order to prove that $I^\gr_{n},~n\in{\mathbb Z}$ is a special inverse system of $N$--K-injective complexes it remains to verify that the complexes 
$I^\gr_{n},~n\in{\mathbb Z}$ are $N$--K--injective, and the short exact sequences
\begin{equation}\label{sh}
0\ra C^\gr_n\ra I^\gr_n\ra I^\gr_{n-1}\ra 0
\end{equation}
are split in each degree, as sequences of $N$--modules.

We shall check that the individual terms of the complexes $I^\gr_{n}$ and $C^\gr_n,~n\in{\mathbb Z}$ are $N$--injective. This implies that the exact sequences (\ref{sh}) are $N$--split in each degree  as exact sequences of $N$--injective modules. 
Moreover, by construction the complexes $I^\gr_{n},~n\in{\mathbb Z}$ are bounded from below, and hence by part 1 of Proposition \ref{sinjprop} $N$--injectivity of the individual terms of these complexes implies $N$--K--injectivity of $I^\gr_{n},~n\in{\mathbb Z}$.

Now we check that the modules  $K^{p,q}$ entering the definition of the complexes $I^\gr_{n},~n\in{\mathbb Z}$ are $N$--injective. Images and cokernels of the differential $d_2$ are analyzed in a similar way. We start with the following simple lemma.
\begin{lemma}\label{sbarm}
Let $M\in\Armb$ be a right $A$-module. Then the bigraded components 
$$
X^{p,q}=\hA(\Baro^{-p}(A,B,A),M)\otimes_A\Baro^q(A,N,A)
$$ 
of the  complex $\sBar(A,N,M)$ may be represented as
$$
X^{p,q}=M\otimes {\rm hom}_{\k}(N\otimes \overline N^{\otimes p},{\k})\otimes_NT^{-q+1}_1(A,N),
$$
where $T^{-q+1}_1(A,N)$ is the quotient of the tensor product 
$$
T^{-q+1}(A,N)=\underbrace{A\otimes_N\ldots \otimes_NA}_{-q+1 ~\mbox{\tiny  times}}
$$ 
by the subspace $\tilde T^{-q+1}(A,N)$,
\begin{eqnarray*}
\tilde T^{-q+1}(A,N)=\{ a_1\otimes\ldots\otimes a_{-q}\otimes a\in T^{-q+1}(A,N) | \qquad \nonumber \\
\qquad \exists s\in\{1,\ldots,-q\}: a_s\in N\}. \nonumber
\end{eqnarray*}
\end{lemma}

\pr \qquad 
First observe that by the definition of the normalized bar resolution \\ $\Bar(A,N,A)$ (see formula (\ref{bar})) the spaces $X^{p,q}$ may be represented as follows: 
\begin{equation}\label{xpq}
X^{p,q}={\rm hom}_B(T_0^{p+1}(A,B),M)\otimes_NT^{-q+1}_1(A,N),
\end{equation}
where $T^{p+1}_0(A,B)$ is the quotient of the tensor product $T^{p+1}(A,B)=\underbrace{A\otimes_B\ldots \otimes_BA}_{p+1 ~\mbox{\tiny  times}}$ by the subspace $\overline T^{p+1}(A,B)$,
$$
\overline T^{p+1}(A,B)=\{ a_0\otimes\ldots\otimes a_p\in T^{p+1}(A,B) |~\exists s\in\{1,\ldots,p\}: a_s\in B\},
$$
and $T^{-q+1}_1(A,N)$ is defined in the formulation of Lemma \ref{sbarm}

Using the triangular decompositions for the algebra $A$ (see condition (vi) in Section \ref{setup}) we can also rewrite the r.h.s. of (\ref{xpq}) as
$$
X^{p,q}={\rm hom}_{\k}(N\otimes \overline N^{\otimes p},M)\otimes_NT^{-q+1}_1(A,N).
$$

Remark that the vector space $N\otimes \overline N^{\otimes p}$ is positively graded and has finite--dimensional graded components while the grading of the vector space $M$ is bounded from above. Therefore by Lemma \ref{tenshom} we have
\begin{equation}\label{xpq1}
X^{p,q}=M \otimes {\rm hom}_{\k}(N\otimes \overline N^{\otimes p},{\k})\otimes_NT^{-q+1}_1(A,N).
\end{equation}

\qed

From Lemma \ref{sbarm} applied to $M=S_A$ it follows that each left $\oppA$--module $K^{p,q}$ may be represented as 
\begin{equation}\label{kpq1}
K^{p,q}=S_A \otimes {\rm hom}_{\k}(N\otimes \overline N^{\otimes p},{\k})\otimes_NT^{-q+1}_1(A,N)\spr M^\prime.
\end{equation}

Now consider the space $K^{p,q}$ as a left $N$--module.
Using formula (\ref{kpq1}) for $K^{p,q}$ and realization (\ref{SA1}) of $S_A$ as a left $N$--module we have the following expression for $K^{p,q}$ as a left $N$--module:
\begin{eqnarray}\label{kpq3}
K^{p,q}=N^*\otimes B \otimes {\rm hom}_{\k}(N\otimes \overline N^{\otimes p},{\k})\otimes_N \qquad \\ \qquad \qquad T^{-q+1}_1(A,N)\spr M^\prime. \nonumber
\end{eqnarray}

Note that the grading of the vector space 
\begin{equation}\label{kpq2}
B \otimes {\rm hom}_{\k}(N\otimes \overline N^{\otimes p},{\k})\otimes_NT^{-q+1}_1(A,N)\spr M^\prime
\end{equation}
is bounded from above since this space is a subspace of
$$
\begin{array}{l}
B \otimes {\rm hom}_{\k}(N\otimes \overline N^{\otimes p},{\k})\otimes_NT^{-q+1}_1(A,N)\otimes_B M^\prime =\\
\\
B \otimes {\rm hom}_{\k}(N\otimes \overline N^{\otimes p},{\k})\otimes \overline B^{\otimes (-q)}\otimes M^\prime 
\end{array}
$$
whose grading is evidently bounded from above.

Applying Lemma \ref{tenshom} to positively graded vector space $N$ with finite--dimensio--nal graded components and the space (\ref{kpq2}), whose grading is bounded from above, we obtain from (\ref{kpq3}) that 
$$
K^{p,q}={\rm hom}_{\k}(N, B \otimes {\rm hom}_{\k}(N\otimes \overline N^{\otimes p},{\k})\otimes_NT^{-q+1}_1(A,N)\spr M^\prime).
$$
The $N$--module in the r.h.s. of the last equality is obviously injective.
This  concludes the proof of $N$--K--injectivity of the complex $\sBaropp (\oppA ,N,M^\prime)$.

(b) The proof of the fact that the complex $\sBaropp (\oppA ,N,M^\prime)$ is relatively to $N$ K--projective is quite similar to the previous one. Namely, we construct a special direct system of relatively to $N$ K--projective complexes converging to $\sBaropp(\oppA,N,M^\prime)$ (see \cite{Sp} or Appendix 2 in this paper for the definition of special direct systems). Then by  Corollary \ref{spdir} the complex $\sBaropp(\oppA,N,M^\prime)$ is relatively to  $N$ K--projective itself.

The required special direct system is defined as follows.   
Consider a set $P^\gr_{n},~n\in {\mathbb Z}$ of complexes in ${\rm Kom}(\oppA)_0$, where $P^\gr_{n}$ is the total complex of the bicomplex $\tau^1_{\leq n}K^{\gr,\gr}$, and $\tau^1_{\leq n}K^{\gr,\gr}$ denotes the truncation from above of the bicomplex  $K^{p,q}$ with respect to the first grading,
$$
\tau^1_{\leq n}K^{\gr,\gr}=\ldots\ra  K^{n-1,\gr}\ra {\rm Ker}~d_1^{n}\ra 0.
$$
Clearly, $P^\gr_{n},~n\in{\mathbb Z}$ is a direct system of complexes converging to $\sBaropp(\oppA,N,M^\prime)$. The proof of the fact that $P^\gr_{n},~n\in{\mathbb Z}$ is a special direct system of relatively to $N$ K--projective complexes is quite similar to the corresponding part of the proof of statement (a) presented in this section.

(c) The cohomology of the complex 
$$
\sBaropp(\oppA,N,M^\prime)=\sBar(A,N,S_A)\spr M^\prime
$$ 
 may be easily calculated with the help of Proposition \ref{SAprop}. 

Indeed, consider $S_A$ as a right $A$--module.  By Proposition \ref{SAprop} $S_A$ is semijective as a right $A$--module, and hence it is a semijective resolution of itself. From part (c) of Proposition \ref{sres} it follows that every semijective resolution of $S_A$ is homotopically equivalent to the zero complex $\ldots \ra 0\ra S_A \ra 0 \ra \ldots$. In particular, the standard semijective resolution $\sBar(A,N,S_A)$ is homotopically equivalent to the 0--complex $\ldots \ra 0\ra S_A \ra 0 \ra \ldots$. This implies that   
\begin{equation}\label{viso}
H^\gr(\sBaropp(\oppA,N,M^\prime))=S_A\spr M^\prime
\end{equation}  
as a vector space. We have to prove that (\ref{viso}) is an isomorphism of left $\oppA$ modules.

In order to do that we explicitly calculate, using spectral sequences, the cohomology of the complex 
$$
\sBaropp(\oppA,N,M^\prime)=\hA(\Baro(A,B,A),S_A)\otimes_A\Bar(A,N,A)\spr M^\prime
$$ 
as the total cohomology of the bicomplex    
$$
K^{p,q}=\hA(\Baro^{-p}(A,B,A),S_A)\otimes_A\Baro^q(A,N,A)\spr M^\prime,~p\geq 0,~q\leq 0.
$$
As before we denote by $d_1$ and $d_2$ the differentials in $K^{p,q}$ induced by the differentials of $\Bar(A,B,A)$ and $\Bar(A,N,A)$, respectively.

The bicomplex $K^{p,q}$ lies in the fourth quadrant of the p--q plane. Therefore the second filtration of this bicomplex defined by
$$
F_{II}^q=\sum_{s\geq q}\sum_{p}K^{p,s}
$$
is regular (see \cite{carteil}, Ch. XV, \S 6). 
The first term of the corresponding spectral sequence may be calculated as the cohomology of the bicomplex $K^{p,q}$ with respect to the differential $d_1$,
$$
E^{p,q}_1=H^p(K^{\gr,q},~d_1).
$$

But for each $q$ the complex $K^{\gr,q}$ may be rewritten, using the definition of the standard normalized bar resolution (see Section \ref{stres}) and the triangular decompositions for the algebra $A$ (see condition (vi) of Section \ref{setup}) , as follows:
\begin{equation}\label{kgrq}
K^{\gr,q}={\rm hom}_N(\Bar(N,{\k},N),S_A)\otimes_NT^{-q+1}_1(A,N)\spr M^\prime,
\end{equation}
where $T^{-q+1}_1(A,N)$ is defined in Lemma \ref{sbarm}.

Now observe that ${\rm hom}_N(\Bar(N,{\k},N),S_A)$ is an $N$--injective resolution of $S_A$ regarded as a right $N$--module. We know that $S_A$ is injective as a right $N$--module (see Proposition \ref{SAprop}), and hence the complex ${\rm hom}_N(\Bar(N,{\k},N),S_A)$ is homotopically equivalent to the 0--complex $\ldots \ra 0\ra S_A \ra 0 \ra \ldots$ as a complex of right $N$--modules. We conclude that
\begin{equation}\label{e1}
\begin{array}{l}
E^{p,q}_1=S_A\otimes_NT^{-q+1}_1(A,N)\spr M^\prime\cdot \delta_{p,0}= \\ 
\\
S_A\otimes_A\Baro^q(A,N,A)\spr M^\prime\cdot \delta_{p,0}
\end{array}
\end{equation}
as a vector space. But from the explicit formula for the differential $d_1$ of the complex (\ref{kgrq}) it follows that (\ref{e1}) is also an isomorphism of left $\oppA$--modules.

Next observe that our spectral sequence degenerates at the second term. By Theorem 5.12, Ch. XV, \S 5 in \cite{carteil} the total cohomology of the bicomplex $K^{p,q}$ is given, as a right $\oppA$--module, by
$$
H^q(\sBaropp(\oppA,N,M^\prime))=E_2^{0,q}=H^q(S_A\otimes_A\Bar(A,N,A)\spr M^\prime).
$$
  
On the other hand from (\ref{viso}) we know that
\begin{equation}\label{viso1}
H^q(\sBaropp(\oppA,N,M^\prime))=S_A\spr M^\prime\cdot \delta_{q,0}
\end{equation}
as a vector space. Now using the explicit form of the differential $d_2$ of the complex $S_A\otimes_A\Bar(A,N,A)\spr M^\prime$ we conclude that (\ref{viso1}) is, in fact, an isomorphism of left $\oppA$--modules. This completes the proof.

\qed


\subsection*{A5. Proof of Theorem \ref{threetor}}

\renewcommand{\thetheorem}{A5.\arabic{theorem}}
\renewcommand{\theequation}{A5.\arabic{equation}}
\renewcommand{\thelemma}{A5.\arabic{lemma}}
\renewcommand{\thedefinition}{A5.\arabic{definition}}

\setcounter{equation}{0}
\setcounter{theorem}{0}
\setcounter{definition}{0}

(1) First we show that the spaces defined by formulas (b) and (c) of Theorem \ref{threetor} are isomorphic. To establish this isomorphism we shall use the standard semijective resolutions constructed in Propositions \ref{res1} and \ref{res2}. More precisely, we shall substitute the standard semijective resolutions $\sBaropp(\oppA,N,M^\prime)$ and $\sBar(A,N,M)$ instead of $S^\gr(M^\prime)$ and $S^\gr(M)$, respectively, into formulas (b) and (c) of Theorem \ref{threetor} and prove that the complexes
$$
M\spr \sBaropp(\oppA,N,M^\prime)
$$
and
$$
\sBar(A,N,M)\spr M^\prime
$$
proposed in Theorem \ref{threetor} for calculation of the space $\stor(M,M^\prime)$ are isomorphic. 

We start with the following simple lemma
\begin{lemma}\label{sbars}
Let $M\in \Almb$ be a right $A$--module. Then the complex $\sBar(A,N,M)$ may be represented as
$$
\sBar(A,N,M)=M\spr \sBar(A,N,S_A).
$$
\end{lemma}

\pr
First observe that by Lemma \ref{sbarm} the bigraded components 
$$
X^{p,q}=\hA(\Baro^{-p}(A,B,A),M)\otimes_A\Baro^q(A,N,A)
$$ 
of the  complex 
$$
\sBar(A,N,M)=\hA(\Baro(A,B,A),M)\otimes_A\Baro(A,N,A)
$$ 
may be represented as
$$
X^{p,q}=M\otimes {\rm hom}_{\k}(N\otimes \overline N^{\otimes p},{\k})\otimes_NT^{-q+1}_1(A,N),
$$
where the spaces $T^{-q+1}_1(A,N)$) are defined in  Lemma \ref{sbarm}. 

Next using Lemma \ref{propspr} we obtain that
$$ 
X^{p,q}=M\spr S_A\otimes {\rm hom}_{\k}(N\otimes \overline N^{\otimes p},{\k})\otimes_NT^{-q+1}_1(A,N).
$$

From  Lemma \ref{sbarm} it follows that  the left $\oppA$--modules 
$$
S_A\otimes {\rm hom}_{\k}(N\otimes \overline N^{\otimes p},{\k})\otimes_NT^{-q+1}_1(A,N)
$$ 
are isomorphic to the bigraded components
$$
\hA(\Baro^{-p}(A,B,A),S_A)\otimes_A\Baro^q(A,N,A)
$$
of the bicomplex $\sBar(A,N,S_A)$. Therefore we have an isomorphism of complexes,
$$
\sBar(A,N,M)=M\spr \sBar(A,N,S_A).
$$
This completes the proof of the lemma.

\qed

Now from the definition of the resolution $\sBaropp(\oppA,N,M^\prime)$, 
$$
\sBaropp(\oppA,N,M^\prime)=\sBar(A,N,S_A)\spr M^\prime,
$$
and Lemma \ref{sbars} we obtain the required isomorphism of complexes,
$$
\begin{array}{l}
\sBar(A,N,M)\spr M^\prime = M\spr \sBar(A,N,S_A)\spr M^\prime = \\
\\
M\spr \sBaropp(\oppA,N,M^\prime).
\end{array}
$$

(2) Next we prove that formulas (a) and (c) of Theorem \ref{threetor} are equivalent. This will complete the proof of Theorem \ref{threetor}. Again we shall use the standard semijective resolutions $\sBaropp(\oppA,N,M^\prime)$ and $\sBar(A,N,M)$.
Substituting these resolutions instead of $S^\gr(M^\prime)$ and $S^\gr(M)$, respectively, into formula (a) of Theorem \ref{threetor} we obtain the following standard complex for calculation of the space $\stor(M,M^\prime)$:
$$
\begin{array}{l}
\sBar(A,N,M)\spr \sBaropp(\oppA,N,M^\prime)= \\
\\
\sBar(A,N,M)\spr \sBar(A,N,S_A)\spr M^\prime
\end{array}
$$

Denote the complex $\sBar(A,N,M)\spr \sBar(A,N,S_A)$ by $Y^\gr$,
$$
Y^\gr=\sBar(A,N,M)\spr \sBar(A,N,S_A).
$$ 
We shall show that $Y^\gr$ is a semijective resolution of $M$. This obviously ensures that formulas (a) and (c) of Theorem \ref{threetor} are equivalent.

The proof of the fact that $Y^\gr$ is a semijective resolution of $M$ is parallel to the proof of Proposition \ref{res1} (see Appendix 3). 

Observe that using the definitions of the standard semijective resolutions (see Propositions \ref{res1} and \ref{res2}) and  Lemma \ref{sbars} the complex $Y^\gr$ may be represented as the total complex of the following four--complex:
\begin{eqnarray}\label{ygrdef}
Y^\gr = \bigoplus_{p\in {\mathbb Z}}\bigoplus_{m+t+v+r=-p}\hA(\Baro^r(A,B,A), \qquad  \qquad \qquad \qquad \qquad  \\ 
 \qquad \qquad \hA(\Baro^t(A,B,A),M)\otimes_A\Baro^{-v}(A,N,A))\otimes_A\Baro^{-m}(A,N,A). \nonumber
\end{eqnarray}
The complex in the r.h.s. of the last equality is a resolution of $M$, i.e. $H^\gr(Y^\gr)=M$, because it consists of standard bar resolutions (see for instance \cite{GM}, Lemma III.7.12). Now by Proposition \ref{resmod} it suffices to verify that $Y^\gr$ is a semijective complex.

First we prove that the complex $Y^\gr$ is K--injective as a complex of $N$--modules. 
Since by part (iii) of Proposition \ref{barprop} the complex $\Baro(A,N,A)$ is homotopically equivalent to $A$ as a complex of $A$--$N$--bimodules, the complex $Y^\gr$ is homotopically equivalent to
\begin{eqnarray*}
\bigoplus_{p\in {\mathbb Z}}\bigoplus_{t+v+r=-p}\hA(\Baro^r(A,B,A),\hA(\Baro^t(A,B,A),M)\otimes_A \\
\qquad \qquad \qquad \qquad \qquad \Baro^{-v}(A,N,A)) \nonumber
\end{eqnarray*}
as a complex of $N$--modules. 

Using the definition of the standard normalized bar resolution (see Section \ref{stres}) and the triangular decompositions for the algebra $A$ (see condition (vi) of Section \ref{setup}) we also obtain that  the last complex, as a complex of $N$--modules, is isomorphic to
\begin{eqnarray}\label{ygr}
\bigoplus_{p\in {\mathbb Z}}\bigoplus_{t+v+r=-p}{\rm hom}_N(\Baro^r(N,{\k},N), \qquad \qquad \qquad \qquad \qquad \qquad \\ 
\qquad \qquad \qquad \hA(\Baro^t(A,B,A),M)\otimes_A\Baro^{-v}(A,N,A)). \nonumber
\end{eqnarray}

Applying again part (iii) of Proposition \ref{barprop} we conclude that the complex (\ref{ygr}) is homotopically equivalent to the complex
$$ 
\bigoplus_{p\in {\mathbb Z}}\bigoplus_{t+r=-p}{\rm hom}_N(\Baro^r(N,{\k},N),\hA(\Baro^t(A,B,A),M)),
$$
that is, in turn, isomorphic to
$$ 
\bigoplus_{p\in {\mathbb Z}}\bigoplus_{t+r=-p}{\rm hom}_N(\Baro^r(N,{\k},N),{\rm hom}_N(\Baro^t(N,{\k},N),M)).
$$

The last complex of is obviously a bounded from below complex of $N$--injective modules. From Proposition \ref{sinjprop} it follows that this complex is also K--injective. This proves that $Y^\gr$ is K--injective as a complex of $N$--modules.

Next we show that $Y^\gr$ is K--projective relative to $N$. By definition we have to prove that for every complex of right $A$--modules $V^\gr\in {\rm Kom}\Armb$, that is homotopically equivalent to zero as a complex of $N$--modules, the complex $\hAd(Y^\gr,V^\gr)$ is acyclic.
Denote this complex by $X^\gr$, 
$$
X^\gr=\HAd(Y^\gr,V^\gr).
$$

We shall apply Lemmas \ref{abovebound} and \ref{belowbound} to calculate the cohomology of this complex. First we  explicitly rewrite the complex $X^\gr$ in components. Using formula (\ref{ygrdef}) for $Y^\gr$ we have:
\begin{equation}\label{be1}
X^\gr=\bigoplus_{q\in {\mathbb Z}}\prod_{l+p=q}\prod_{m+t+v+r=p}X^{l,r,t,v,m},
\end{equation} 
where 
\begin{eqnarray*}
X^{l,r,t,v,m}=\hA(\hA(\Baro^r(A,B,A),\hA(\Baro^t(A,B,A),M)\otimes_A \qquad \\
\qquad \qquad \qquad \Baro^{-v}(A,N,A))\otimes_A\Baro^{-m}(A,N,A),V^l).
\end{eqnarray*}

Note that $X^\gr$ is an ``almost'' five--complex with components $X^{l,r,t,v,m}$. 
We denote by $d_1$ the differential in $X^\gr$ induced by the differential of the complex $V^\gr$, and by $d_2,~d_3$ and $d_4,~d_5$ the differentials in $X^\gr$ induced by the differentials of the complexes $\Bar(A,B,A)$ and $\Bar(A,N,A)$, respectively. The total differential $d$  of $X^\gr$ is the sum $d=d_1+d_2+d_3+d_4+d_5$.

First using Lemma \ref{abovebound} we calculate the cohomology of the complex $X^\gr$ with respect to the differential $d'=d_1+d_2+d_3$. In order to do that we note that by the definition of the normalized bar resolution $\Bar(A,N,A)$ (see formula (\ref{bar})) the complex (\ref{be1}), with the differentials $d_4,~d_5$ forgotten, may be rewritten as 
\begin{eqnarray*}
X^\gr=\hA(\hA(\Baro^r(A,B,A),\hA(\Baro^t(A,B,A),M)\otimes_A \qquad \\
\qquad \qquad \qquad \Baro^{-v}(A,N,A))\otimes_AT^{m+1}_0(A,N),V^l),
\end{eqnarray*}
where the spaces $T^{n+1}_0(A,N)$ are defined in the proof of Proposition \ref{res1} (see formula \ref{star} in Appendix 3).

Now recall that the complex $V^\gr$ is homotopically equivalent to zero as a complex of $N$--modules. Observe also that the functor ${\rm hom}_N$ preserves this property. Therefore for every $r,t,v,m\in {\mathbb Z}$ the complex $X^{\gr,r,t,v,m}$ is homotopic to zero. Note also that $X^{\gr,r,t,v,m}=0$ for $r,t>0$, and hence from Lemma \ref{abovebound}, applied to the complex $X^\gr$ equipped with the differential $d'=d_1+d_2+d_3$, it follows that $H^{\gr}(X^\gr,d')=0$.

Next, we apply Lemma \ref{belowbound} to the complex $(X^{\gr},d)$, the differential $d$ being split as follows $d=d'+d_4+d_5$. It suffices to remark that $X^{l,r,t,v,m}=0$ for $v,m<0$, and $H^{\gr}(X^\gr,d')=0$ as we proved above. Therefore by  Lemma \ref{belowbound} $H^\gr(X^\gr,d)=0$. This proves that the complex $Y^\gr$ is K--projective relative to $N$. But we have already proved that this complex is K--injective as a complex of $N$--modules. We conclude that $Y^\gr$ is a semijective complex. This completes the proof of Theorem \ref{threetor}

\qed

\end{document}